%% file: main.tex
\documentclass[review,hidelinks,onefignum,onetabnum]{siamart250211}

\input{ex_shared}

\ifpdf
\hypersetup{
  pdftitle={Exact operator inference with minimal data}{H. Rosenberger, B. Sanderse, and G. Stabile},
  pdfauthor={H. Rosenberger, B. Sanderse, and G. Stabile}
}
\fi

\input{preamble/packages}

\input{preamble/macros}

\input{preamble/rank_suff_macros}

\input{preamble/main_only_macros}
\graphicspath{{./figures/}{./figures/inkscape/}}

\setuptodonotes{fancyline}

\input{revision/report_commands}

\begin{document}

\maketitle

\begin{abstract}
This work introduces a novel method to generate snapshot data for operator inference that 
guarantees
the exact  
reconstruction of intrusive projection-based reduced-order models (ROMs).
To ensure exact reconstruction, the operator inference least squares matrix must have full rank, without regularization.
Existing works have achieved this full rank using
heuristic strategies to generate snapshot data and a-posteriori checks on full rank, 
but without a guarantee of success.
Our novel snapshot data generation method provides this guarantee
thanks to
two key ingredients:
first we identify ROM states that induce full rank, then we generate snapshots corresponding to exactly these states by simulating multiple trajectories for only a single time step.
This way, the number of required snapshots is minimal and orders of magnitude lower than typically reported with existing methods.
The method avoids non-Markovian terms and does not require \reprojection.
Since the number of snapshots is minimal, the least squares problem simplifies to a linear system that is numerically more stable.
In addition, because the inferred operators are exact, properties
of the intrusive ROM operators 
such as symmetry or skew-symmetry are preserved.
Numerical results for differential equations involving 2nd, 3rd and 8th order polynomials demonstrate that the novel 
snapshot data generation 
method
leads to
exact reconstruction of the intrusive 
reduced order models.

\end{abstract}

\begin{keywords}
operator inference, data-driven modeling, \nonintrusive\ model reduction, reduced basis method
\end{keywords}

\begin{MSCcodes}
65Y99, 65F22, 35R30, 65D05
\end{MSCcodes}


\tableofcontents

\input{sections/introduction}

\input{sections/preliminaries}

\section{Exact operator inference}
\label{sec:rank suff snapshot data generation}

\newcommand{\cx}{\bar\xv}
\newcommand{\cu}{\bar\uv}

In this work, we propose a novel method to generate snapshots that is a priori guaranteed to result in a full rank 
matrix $\Pm$. Our method has two key ingredients. First, in Section \ref{sec:basis of the feature space}, we leverage the 
knowledge about how the matrix $\Pm$  depends on the ROM states $\bx_k$ and input signals $\uv_k$ via \eqref{eq:def p matrix and dot breve x matrix} and \eqref{eq:def p vec}
to find pairs of ROM states $\cx_\ind\in\mbR^n$ and input signals $\cu_\ind\in\mbR^\Nu$, $\ind=1,\dots,\nf$ that ensure full rank of $\Pm$. Second, in Section \ref{sec:single step ensemble}, we generate
time derivative snapshots $\dot\cx_\ind$ corresponding to
these rank-ensuring pairs $(\cx_\ind,\cu_\ind)$ by simulating $\nf$ trajectories for only a single time step each, with initial conditions $\cx_\ind$ and input signals $\cu_\ind$.

\subsection{Rank-ensuring 
ROM states
and input signals}
\label{sec:basis of the feature space}

The motivation for our choice of rank-ensuring ROM states is to isolate single entries of the aggregated polynomial $\pv(\tilde \xv,\uv)$ defined in \eqref{eq:def p vec}. If we could find inputs $\tilde \xv$ and $\uv$ such that $\pv(\tilde \xv,\uv)$ is 1 in the $j$-th entry and zero in all other entries for all $j=1,\dots,\nf$, then the product $\hat \Om \pv(\tilde \xv,\uv)$ would equal the $j$-th column of the operator $\hat \Om$. Hence, we could easily infer all columns of the operator $\hat \Om$.
\R{2 example 1}
\revtwo{
\begin{example}
\label{example 1}
  If $\inds=\{1\}$ and $\Nu=0$, so $\pv(\tilde \xv) = \revfour{\tilde \xv}$, we can achieve such an isolation of entries by using the unit vectors in $\mbR^n$, here denoted by $\onehot_j$, $j=1,\dots,n$, as inputs $\tilde \xv$.  
\end{example}}

\R{2 example 2}
\revtwo{
\begin{example}
\label{example 2}
    If $\inds=\{2\}$ and $\Nu=0$, so $\pv(\tilde \xv) = \tilde \xv^2$, perfect isolation of entries is not possible. For example, when choosing $\tilde \xv \;\revone{= [\tilde x_1,\dots,\tilde x_n]^T} = \onehot_i + \onehot_j$ such that $\pv\revfour{(\tilde \xv)}$ is 1 in the entry $\revfour{\tilde x_i \tilde x_j}$ with $i\neq j$, also the entries $\revfour{\tilde x_i^2}$ and $\revfour{\tilde x_j^2}$ are 1. However, the columns of $\hat \Om$ corresponding to $\revfour{\tilde x_i \tilde x_j}$, $\revfour{\tilde x_i^2}$ and $\revfour{\tilde x_j^2}$ can still be inferred when employing the data $\hat\Om\pv(\onehot_i+\onehot_j)$, $\hat\Om\pv(2\onehot_i)$ and $\hat\Om\pv(2\onehot_j)$.
\end{example}
}

Based on this motivation, 
we define
for arbitrary $i\in\mbN^+$
the set $\setXi$ that consists of all sums of $i$ of these vectors,
\begin{align}
    \setXi := 
    \left\{
       \sum_{j=1}^i \vt a_j \Bigg| \vt a_1,\dots,\vt a_i \in\{\onehot_1,\dots,\onehot_n\}
    \right\}
    \h \text{ and consistently } \setXx 0 := \{\vt 0\in\mbR^n\}
    .
    \label{eq:def set X i}
\end{align}
The main result of this article is: the vectors in $\cup_{i\in\inds} \setXi$  induce full rank of the matrix $\Pm$ as defined in \eqref{eq:def p matrix and dot breve x matrix} and \eqref{eq:def p vec}, for any \degreeset\ $\inds\subset\mbN_0$ and $\Nu = 0$.

If $\Nu > 0$, we achieve full rank of $\Pm$ 
by pairing these ROM states $\cx_j$ with zero input signals, $(\cx_j,\vt 0)$, and adding
the pairs $(\vt 0, \uv_j)$ with zero ROM state and $\uv_j$, $j=1,\dots,\Nu$ the unit vectors in $\mbR^\Nu$.
This way we get the rank-ensuring snapshot pairs:
\begin{align}
    \left\{(\cx_j,\cu_j)\right\}_{j=1}^\nf = \left\{(\cx_j,\vt 0) | \tx_j\in\cup_{i\in\inds} \setXi\right\} \cup \left\{(\vt 0, \uv_j)\right\}_{j=1}^\Nu.
    \label{eq:initial condition and input signal pairs}
\end{align}

\R{2 example 3}
\revtwo{
\begin{example}
\label{example 3}
For
$n=2$, $\Nu = 2$ and $\inds = \{1,2\}$ (such that $\np = 5$, $\nf = 7$), we get the pairs
\begin{align*}
    \underbrace{
    \left(\begin{bmatrix} 1 \\ 0 \end{bmatrix}, \begin{bmatrix} 0 \\ 0 \end{bmatrix}\right),
    \left(\begin{bmatrix} 0 \\ 1 \end{bmatrix}, \begin{bmatrix} 0 \\ 0 \end{bmatrix}\right)
    }_{\left\{(\cx_j,\vt 0) | \tx_j\in \mathcal X^1\right\}},
    \underbrace{
    \left(\begin{bmatrix} 2 \\ 0 \end{bmatrix}, \begin{bmatrix} 0 \\ 0 \end{bmatrix}\right),
    \left(\begin{bmatrix} 1 \\ 1 \end{bmatrix}, \begin{bmatrix} 0 \\ 0 \end{bmatrix}\right),
    \left(\begin{bmatrix} 0 \\ 2 \end{bmatrix},
    \begin{bmatrix} 0 \\ 0 \end{bmatrix}\right)
    }_{\left\{(\cx_j,\vt 0) | \tx_j\in \mathcal X^2\right\}},
    \underbrace{
    \left(\begin{bmatrix} 0 \\ 0 \end{bmatrix}, \begin{bmatrix} 1 \\ 0 \end{bmatrix}\right),
    \left(\begin{bmatrix} 0 \\ 0 \end{bmatrix}, \begin{bmatrix} 0 \\ 1 \end{bmatrix}\right)
    }_{\left\{(\vt 0, \uv_j)\right\}_{j=1}^\Nu},
\end{align*}
that result in the full-rank matrix
\begin{align*}
    \Pm = \begin{bmatrix} 
    1 & 0 & 2 & 1 & 0 & 0 & 0 \\
    0 & 1 & 0 & 1 & 2 & 0 & 0 \\
    1 & 0 & 4 & 1 & 0 & 0 & 0 \\
    0 & 0 & 0 & 1 & 0 & 0 & 0 \\
    0 & 1 & 0 & 1 & 4 & 0 & 0 \\
    0 & 0 & 0 & 0 & 0 & 1 & 0 \\
    0 & 0 & 0 & 0 & 0 & 0 & 1 \\
    \end{bmatrix}  \h \text{ where } \h 
    \pv(\tx,\uv) = \revfour{\begin{bmatrix} \tilde x_1 \\ \tilde x_2 \\ \tilde x_1^2 \\ \tilde x_1\tilde x_2 \\ \tilde x_2^2 \\ u_1 \\ u_2 \end{bmatrix}}.
\end{align*}
\end{example}
}

\R{2:full rank theo}
\revtwo{
In general, the following theorem holds.
\begin{mytheo}[Full rank of matrix $\Pm$]
\label{theo:full rank P}
For any $n\in\mbN^+$, $\Nu\in\mbN_0$ and $\inds\subset\mbN_0$,
the matrix
$\Pm := 
    \left[ 
    \pv(\bx_0,\uv_0) \;\;\dots\;\;
    \pv(\bx_{K-1},\uv_{K-1})
    \right]
    \; \in\mbR^{\nf\times K}$ 
defined in \eqref{eq:def p matrix and dot breve x matrix} with $K=\nf$, $\pv(\tx,\uv) = \begin{bmatrix}
        \left[\tx^i\right]_{i\in\inds} \\ \uv
    \end{bmatrix}$ as defined in \eqref{eq:def p vec} 
and with the ROM states and input signals 
$
    \left\{(\cx_j,\cu_j)\right\}_{j=1}^\nf = \left\{(\cx_j,\vt 0) | \tx_j\in\cup_{i\in\inds} \setXi\right\} \cup \left\{(\vt 0, \uv_j)\right\}_{j=1}^\Nu
    $
defined in \eqref{eq:initial condition and input signal pairs}
has full rank.
\end{mytheo}
The proof is given in Appendix \ref{sec:full rank proof}.
}


Note that this choice of \revtwo{ROM states and input signals} results in a matrix $\Pm$ that has as many columns as rows. This means that our choice uses the minimal number of snapshots necessary to 
achieve
full rank.
In addition, 
because $\Pm$ is square,
the operator inference minimization problem \eqref{eq:opinf minimi matrices} is equivalent to the matrix-valued linear system
\begin{align}
    \hat \Om \Pm = \dot \bX.
    \label{eq:opinf linear system}
\end{align}
Since the matrix $\Pm$ has a smaller condition number than the matrix $\Pm\Pm^T$ of the normal equations 
\begin{align}
    \hat \Om \Pm\Pm^T = \dot\bX \Pm^T
\end{align}
corresponding to the least squares problem \eqref{eq:opinf minimi matrices}, \R{2: this linear system}\revtwo{the linear system \eqref{eq:opinf linear system}} is in general
less sensitive to error amplifications.

\subsection{Snapshot data generation via single-step ensemble}
\label{sec:single step ensemble}

To formulate the operator inference minimization problem \eqref{eq:opinf minimi vecs}, we not only need the pairs of ROM states and input signals, denoted by $(\cx_j,\cu_j)$, $j=1,\dots,\nf$ defined in 
\eqref{eq:initial condition and input signal pairs}, but also approximations of the ROM state time derivative $\dot \cx_j$ that correspond to these pairs.
In 
existing works on operator inference \cite{PEHERSTORFER2016data,peherstorfer2020sampling,mcquarrie2021data,goyal2023guaranteed}, 
triples of these quantities are
obtained by projection of the corresponding FOM vector that are
generated along trajectories of the FOM. 
Hence, the FOM state $\xv_j$ at a given time step $t_j$ depends on the FOM state $\xv_{j-1}$ at the previous time step $t_{j-1}$ and the FOM dynamics \eqref{eq:fom dynamics}. These FOM dynamics are unknown (otherwise we would not need to perform operator inference). Consequently, we cannot predict or control which FOM states are attained along such trajectories. In particular, we cannot guarantee that FOM states corresponding to the ROM states $\cx_j$ described in the previous section are attained.

In our approach,
we do not generate FOM snapshot data
along one (or few) long trajectories but instead use an ensemble of many
short trajectories, each containing \textit{only a single time step}. 
Each trajectory is initialized with 
\begin{align}
    \text{ the 
    rank-ensuring
    initial condition $\x_0^\ind = \Vn \cx_\ind$, and input signal $\uv_0^\ind = \cu_\ind$, $\ind=1,\dots,\nf$, }
    \label{eq:def well-chosen initials and inputs}
\end{align}
\revtwo{with $\cx_\ind$ and $\cu_\ind$ as defined in \eqref{eq:initial condition and input signal pairs},}
and integrated in time for only a single explicit
Euler step of size $\Delta t$. This approach results in FOM snapshots
\begin{align}
    \x_1^\ind = \x_0^\ind + \Delta t \f(\x_0^\ind,\uv_0^\ind) \h\h\h \ind=1,\dots,\nf.
    \label{eq:single time step snapshot}
\end{align}
From each pair of FOM snapshots $\x_0^\ind$ and $\x_1^\ind$, we can compute the ROM time derivative snapshot
\begin{align}
    \dot \cx_\ind = 
    \Vn^T \frac{\x_1^\ind - \x_0^\ind}
    {\Delta t} \h\h\h \ind=1,\dots,\nf.
\label{eq:rank suff time derivative}
\end{align}
The resulting
triples \R{4:consistent triple 2}$(\revfour{\cx_\ind,\dot\cx_\ind},\uv_\ind)_{\ind = 1}^\nf$
form the input
to the operator inference optimization problem \eqref{eq:opinf minimi vecs}.
Since the matrix $\Pm$ \eqref{eq:def p matrix and dot breve x matrix} only depends on the snapshots 
$\cx_\ind$ and $\uv_\ind$, $\ind = 1,\dots,\nf$, it does not depend
on the unknown FOM operators. Hence, we can 
guarantee full rank of $\Pm$ by 
the choice of the initial conditions and input signals.

Note 
that our novel approach for generating snapshot data for operator inference does not affect the snapshot data used for POD, so the POD basis $\Vn$ remains unchanged. \R{2:improve staccato 2}\revtwo{For example, the implicit Euler method is used to generate POD snapshot data in Section \ref{sec:ice sheet model}.}

The idea to generate snapshot data in multiple trajectories, each containing only a single time step, has already been used in numerical experiments in \cite{PEHERSTORFER2016data}. However, randomized instead of well-chosen initial conditions 
are
used in that work. \R{1:mention uy}\revone{Later, FOM snapshots projected onto the ROM space have been used as initial conditions in context of lifted nonlinearities \cite{qian2020lift} and noisy data \cite{uy2023active}.
\R{3:a priori}\revthree{
All three approaches enable to check the full rank a priori, but do not guarantee full rank.
}
Indeed, far more than $\nf$ trajectories have been used to achieve full rank of $\Pm$ in the numerical experiments in these works. 
}

The complete method is summarized in 
Algorithm \ref{alg:exact opinf}. 

\begin{algorithm}
\caption{Exact operator inference}
\label{alg:exact opinf}

\begin{algorithmic}[1]
\Procedure{exactOpInf}
{
\newline\phantom{---} ROM basis $\Vn\in\mbR^{N\times n}$,
\newline\phantom{---} \degreeset\ $\inds\subset \mbN_0$,
\newline\phantom{---} input signal dimension $\Nu\in\mbN_0$
\newline\phantom{---} time step size $\Delta t$
    \newline}
    \LineComment{Generate 
    rank-ensuring
    snapshot data.}
    \State Initialize $\Pm\gets\emptyset$
    \State Initialize $\dot \bX \gets \emptyset$
    \revtwo{
    \For{$i\in\inds$}
        \LineComment{Generate ROM state sets $\setXi$ as defined in 
        \eqref{eq:def set X i} 
        (see Appendix \ref{appendix: algo gen rom states})}
        \State $\setXx i\gets$ 
        \Call{GenerateROMstateSet}{ROM dimension $n$, polynomial degree $i$} 
        \label{alg line:gen rom state set}
    \EndFor}
    \For{$\cx\in \cup_{i\in\inds} \setXi$ 
    }
    \label{alg line: set X i}
        \State $\xv_1 \gets$ \textproc{FOM}(
initial condition $\xv_0 = \Vn\cx$,
\newline\phantom{--------------------------} input signal $\uv = \vt 0$,
\newline\phantom{--------------------------} time step size $\Delta t$, end time $T = \Delta t$, explicit Euler time discretization
)
    \State $\Pm \gets [\Pm, \;\; \pv(\tx_0,\vt 0)]$
    \State $\dot \bX \gets [\dot\bX, \;\; \frac{\Vn^T\xv_1 - \cx}{\Delta t}]$
    \EndFor

    \For{$j=1,\dots,\Nu$}
        \State $\uv \gets$ $j$-th unit vector in $\mbR^\Nu$
        \State $\xv_1 \gets$ \textproc{FOM}(
initial condition $\xv_0 = \vt 0$,
\newline\phantom{--------------------------} input signal $\uv$,
\newline\phantom{--------------------------} time step size $\Delta t$, end time $T = \Delta t$, explicit Euler time discretization
)
    \State $\Pm \gets [\Pm, \;\; \pv(\vt 0,\uv)]$
    \State $\dot \bX \gets [\dot\bX, \;\; \frac{\Vn^T\xv_1 - \vt 0}{\Delta t}]$
    \EndFor

    \State Solve linear system \eqref{eq:opinf linear system}
    \State
    \revfour{\Return $\hat \Om$}
    \label{alg line: return hat o}
\EndProcedure
\end{algorithmic}

\end{algorithm}

\subsection{Exact reconstruction of the intrusive ROM operators}

\revtwo{
By combining the rank-ensuring initial conditions and input signals in Section \ref{sec:basis of the feature space} and the single-step ensemble in Section \ref{sec:single step ensemble}, our novel approach exactly reconstructs the intrusive ROM operators, as stated by the following theorem.
\R{2:theorem exact reconstruction}
\begin{mytheo}[Exact reconstruction of intrusive ROM operators]
\label{theo:exact reconstruction}
    For any $n\in\mbN^+$, $\Nu\in\mbN_0$, $\inds\subset\mbN_0$ and $\Delta t>0$,
    the 
    operator inference minimization problem \eqref{eq:opinf minimi vecs} with snapshot data $(\bx_k,\dot\bx_k,\uv_k)=(\cx_\ind,\dot\cx_\ind,\cu_\ind)$ defined in \eqref{eq:initial condition and input signal pairs} and \eqref{eq:rank suff time derivative} has a unique solution. This solution is equivalent to the intrusive ROM operators defined in \eqref{eq:intrusive ops1} and \eqref{eq:intrusive ops2}.
\end{mytheo}
\begin{proof}
We first show that the FOM snapshot data generated in our proposed approach is equivalent to snapshot data generated from the intrusive ROM itself, which constitutes the ideal inference data.
For this purpose, we 
insert the definition of the 
rank-ensuring
initial conditions 
$\cx_\ind$
\eqref{eq:def well-chosen initials and inputs} 
and the explicit Euler FOM time step \eqref{eq:single time step snapshot}
into the definition of the explicit Euler time derivative snapshot 
\eqref{eq:rank suff time derivative}
to find
\begin{align}
    \dot \cx_\ind = 
        \frac 1 {\Delta t}\Vn^T\left( \x_0^\ind + \Delta t \f(\x_0^\ind,\uv_0^\ind) -\x_0^\ind\right)
    = \Vn^T \f(\x_0^\ind,\uv_0^\ind)
    = \Vn^T \f(\Vn \cx_\ind,\uv_0^\ind) = \tilde\Om \pv(\cx_\ind,\uv_0^\ind).
    \label{eq:exact rhs reconstruction}
\end{align}
Hence, \eqref{eq:opinf minimi vecs} simplifies to
\begin{align}
        \min_{\hat \Om} \sum_{\ind = 1}^\nf \|\hat \Om \pv(\cx_\ind,\uv_0^\ind) - \tilde \Om \pv(\cx_\ind,\uv_0^\ind)\|_2^2,
\end{align}
or equivalently
\begin{align}
    \min_{\hat \Om} \|\hat \Om \Pm - \tilde \Om \Pm\|_F^2.
\end{align}
Since $\Pm$ has full rank as shown in Theorem \ref{theo:full rank P}, 
the solution $\hat\Om = \tilde \Om$ is unique, so the method exactly reconstructs the intrusive ROM operator.
\end{proof}
}

Note that
our approach not only guarantees the 
full rank
of $\Pm$ but also eliminates two other error sources common in operator inference. Firstly, the consistency of using explicit Euler schemes, both for the snapshot generation \eqref{eq:single time step snapshot} and the computation of the time derivative snapshots 
\eqref{eq:rank suff time derivative} avoids the introduction of time discretization 
inconsistencies.
Secondly, by performing only a single explicit Euler time step on each trajectory we avoid 
introducing non-Markovian terms that would cause the projected FOM snapshot data \eqref{eq:proj fom snapshots} to deviate from actual ROM snapshot data as described in \cite{peherstorfer2020sampling}.

Since our approach exactly reconstructs intrusive ROM operators, it also preserves any structure present in these intrusive operators such as symmetry and skew-symmetry \cite{sanderse2020non}.
In contrast to other works \cite{Koike_Qian_2024,gruber2023canonical,goyal2023guaranteed},
these structures are preserved automatically
without the need to enforce their preservation explicitly.

\subsection{Assumptions on the FOM solver}
\label{sec:assumptions on fom solver}

Operator inference is a \nonintrusive\ model reduction concept developed for FOM solvers that do not provide the access to FOM operators required by intrusive model reduction. Consequently, imposing strong assumptions on the FOM solver contradicts the core motivation for operator inference.
In our proposed approach, we only impose two assumptions.

First, we
assume that specific initial conditions can be imposed; concretely, linear combinations of the ROM basis modes
\eqref{eq:def well-chosen initials and inputs}. 
\R{1:physically realizable}\revone{For POD bases, these modes are linear combinations of FOM snapshots. Since these snapshots are physically realizable, the specific initial conditions are physically realizable as well in most cases. Exceptions include models that only allow a certain range of state values, for example positivity of temperature, or require specific boundary conditions. These exceptions are beyond the scope of this article since they 
typically require some additional treatment already on the intrusive ROM level, for example 
the introduction of a lifting function \cite{rosenberger2023no} to satisfy boundary conditions.}

Second, we assume that explicit Euler can be chosen as time discretization method to generate the operator inference snapshot data.
\R{1:higher-order rk}
\revone{
However, this restriction can be relaxed to higher-order Runge-Kutta schemes if input signals are interpolated between time steps.
These schemes can be written as explicit Euler schemes whose right-hand side has higher polynomial degree \cite{peherstorfer2020sampling}.
Consequently, our approach can be applied to explicit Runge-Kutta methods (and multi-step methods that use an explicit Runge-Kutta method to compute the first time step).
}
For implicit methods, on the other hand, the full rank of $\Pm$ might not be guaranteed. For example, if the FOM time steps \eqref{eq:single time step snapshot} are replaced by implicit Euler steps,
\begin{align}
    \xv_1^\ind = \xv_0^\ind + \Delta t \fv(\xv_1^\ind,\uv_1^\ind),
\end{align}
then the matrix $\Pm$ should, for time discretization consistency, be composed of vectors $\pv(\xv_1^\ind,\uv_1^\ind)$ instead of $\pv(\xv_0^\ind,\uv_0^\ind)$. However, in contrast to the initial conditions $\xv_0^\ind$, we do not have full control over the states $\xv_1^\ind$ after one time step. Consequently, we generally cannot guarantee the full rank of $\Pm$.
\R{1:alternative time disc error}
\revone{Alternatively, we can construct the matrix $\Pm$ from the vectors $\pv(\xv_0^\ind,\uv_0^\ind)$. This way, we guarantee the full rank of $\Pm$ but introduce a time discretization inconsistency that impedes exact reconstruction of the intrusive ROM operators.}

Note that this requirement on the time discretization only affects the generation of operator inference snapshot data, but neither the generation of POD snapshot data, nor the simulations that are performed with the resulting ROM.
\R{2:improve staccato}\revtwo{
The POD snapshot data can be generated from any FOM time discretization, and -
like intrusive ROM operators - the inferred operators can be used for any ROM time discretization.}

\R{1:delta t as input}
\revone{Note further that the time step size $\Delta t$ is provided as an input in Algorithm \ref{alg:exact opinf} only for ease of presentation. Our proposed approach is also applicable to FOM solver that do not allow the users to specify the time step size themselves. In such cases, $\Delta t$ is not provided as an input and the end time $T$ cannot be chosen to be equal to $\Delta t$. Consequently, we can generally not guarantee that the FOM solver only performs a single time step, but by setting $T$ to a small value, we can limit the number of excess time steps. From the solver output, we then only use the state snapshot $\x_1$ after the first time step and the corresponding time $t_1$ to compute the time step size $\Delta t$.}

\R{1:reproj alike}
\revone{While our approach differs decisively from the \reprojection\ approach \cite{peherstorfer2020sampling} in terms of control over which snapshots are generated, both approaches are very similar in terms of intrusiveness. In fact, our approach can be seen as the \reprojection\ approach performed for $\nf$ separate trajectories, but stopped after a single time step and before the first \reprojection\ step. Consequently, the described assumptions on the FOM solver are almost the same as for the \reprojection\ approach. The only exception are implicit time discretization schemes for linear models which are allowed in \reprojection\ but not in our approach because they inhibit control over the generated snapshots.}

\subsection{Remarks on changing ROM dimension}
\label{sec:rem changing rom dim}
 
\newcommand{\nsmall}{{n_-}}
\newcommand{\nbig}{{n_+}}

When ROMs are used in practice, it can be useful to change their dimension.
For instance, the
ROM dimension 
might
be increased to 
improve accuracy or reduced to 
improve computational 
efficiency.
The proposed data generation method is flexible towards such dimension changes. Concretely, the set $\mat X_n$ of FOM snapshots $\xv_1^s$ as defined in \eqref{eq:single time step snapshot} for any ROM dimension $n\in\mbN$ is a subset of the corresponding set $\mat X_{n+1}$ of FOM snapshots for dimension $n+1$. Generally, $\mat X_\nsmall\subset \mat X_n \subset \mat X_\nbig$ for all $\nsmall < n < \nbig$. 
Consequently, we do not need to generate new data to infer ROM operators of dimension $\nsmall$ when the snapshot data $\mat X_n$ is already generated,
but
can select the subset $\mat X_\nsmall$ of $\mat X_n$.
Likewise, we do not need to generate the complete snapshot data $\mat X_\nbig$ when $\mat X_n$ is given,
but
can completely reuse $\mat X_n$ and only have to generate $\mat X_\nbig \setminus \mat X_n$ additionally.

Note that \reprojection\ data \cite{peherstorfer2020sampling} does not have this flexibility because the \reprojection\ process is specific to a fixed ROM dimension.

    \subsection{Choice of time step size $\Delta t$}
    \label{sec:choice of delta t}
    \revtwo{
    As stated in Theorem \ref{theo:exact reconstruction}, our proposed method achieves exact reconstruction for any $\Delta t>0$.
    Indeed,
    in exact arithmetics, the size of $\Delta t$ also does not affect the accuracy of the inferred operators since the effect of $\Delta t$ in computing $\xv_1^\ind$ in \eqref{eq:single time step snapshot} cancels out when computing the time derivative snapshot $\dot \bx_1^\ind$ in 
    \eqref{eq:rank suff time derivative}.
    }
    In floating point arithmetics, however, the accuracy of the inferred operators depends on how accurate the intermediate results can be represented as floating-point numbers.
    \R{1:delta t by solver}
    \revone{If $\Delta t$ is determined by the FOM solver, we have to accept the inaccuracy in the ROM operator reconstruction induced by this limitations of floating-point number representation. However, if $\Delta t$ can be specified by the user, we can choose $\Delta t$ such that these limitations are mitigated.}

    \revone{
    Concretely, the differences $\bx_1^\ind - \bx_0^\ind$ cannot be represented well as floating-point numbers if $\|\tilde\Om\|_2$ 
    is very large or very small.}
    To counteract this effect, $\Delta t = \frac 1 {\|\tilde\Om\|_2}$ would be a good choice. Since we do not know $\tilde\Om$ a priori, we can only approximate $\|\tilde\Om\|_2$ based on the POD snapshot data.
 
    In this work, we use the estimate
    \begin{align}
        \Delta t = \left(\max_{k=0,\dots,K_{POD}-1} \frac{\left\|\vv_1^T\frac{x_{k+1}-x_k}{t_{k+1}-t_k}\right\|_2}{\left\| \pv(\vv_1^T x_k,u(t_k))\right\|_2}\right)^{-1},
        \label{eq:dt estimate}
    \end{align}
    based on the lower bound
    \begin{align}
        \|\tilde\Om\|_2 := \sup_{\vt y\in\mbR^\nf} \frac{\|\tilde\Om\vt y\|_2}{\|\vt y\|_2} \geq \sup_{\tx\in\mbR^n,\, \uv\in\mbR^\Nu} \frac{\|\tilde\Om\pv(\tx,\uv)\|_2}{\|\pv(\tx,\uv)\|_2} \geq \max_{k=0,\dots,K_{POD}-1} \frac{\|\tilde\Om\pv(\Vn^T\x_k,\uv(t_k))\|_2}{\|\pv(\Vn^T\x_k,\uv(t_k))\|_2},
        \label{eq:choice of time step 1}
    \end{align}
    and the approximation
    \begin{align}
        \tilde \Om \pv(\Vn^T\x_k,\uv(t_k)) \approx 
        \Vn^T\frac{\x_{k+1}-\x_k}{t_{k+1}-t_k}, \h k = 0,\dots,K_{POD}-1,
        \label{eq:choice of time step 2}
    \end{align}
    for ROM dimension $n=1$\R{1 v_1}\revone{, so $\Vn = \vt v_1$}.
    
\revtwo{Note that we do not need to take stability into account when determining the time step size $\Delta t$ since we only perform a single time step per trajectory.}

\section{Numerical experiments}
\label{sec:num exps}

The following numerical experiments demonstrate that Algorithm \ref{alg:exact opinf} exactly reconstructs the corresponding intrusive ROM operators.
Many existing works on operator inference assess the accuracy of the inferred ROMs based on trajectory errors. 
As shown in \cite{sawant2023physics}, trajectory errors can be misleading
because they are sensitive to the choice of trajectories.
Furthermore, the objective of operator inference is to reconstruct the intrusive ROM, so this intrusive ROM should serve as ground truth, not the underlying FOM
\R{2:fom errors}
\revtwo{(see Appendix \ref{appendix:rom state errors})}.
We therefore assess the reconstruction accuracy with relative operator errors $\frac{\|\hat\Om-\tilde\Om\|_F}{\|\tilde \Om\|_F}$.

\newcommand{\xs}{x}
\newcommand{\xis}{\xi}
\newcommand{\uus}{u}

The code is available on \url{https://github.com/h3rror/exactOpInf}.

\begin{figure}
\centering
\begin{minipage}[c]{0.49\linewidth}
\centering
\includegraphics[width=\textwidth]{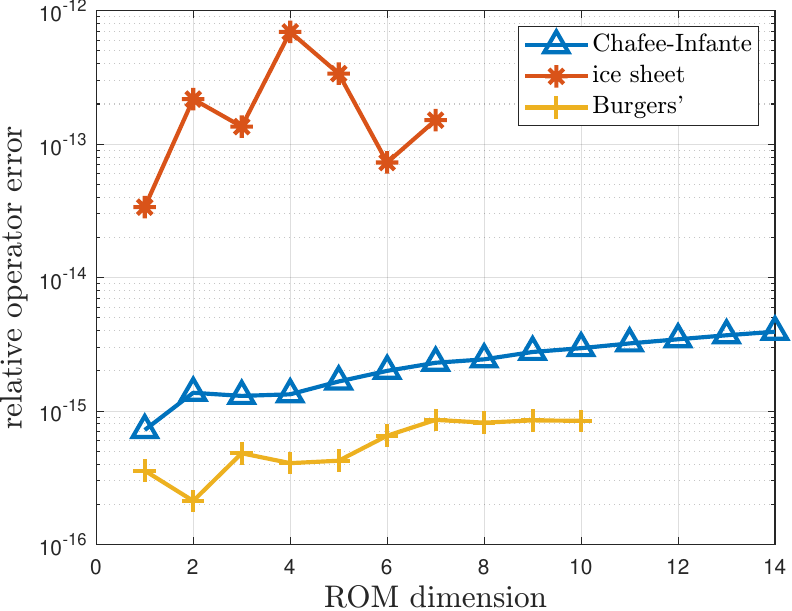}
\caption{
\revfour{
Relative operator errors $\frac{\|\hat \Om - \tilde \Om\|_F}{\|\tilde \Om\|_F}$.}
}
\label{fig:operator errors}
\end{minipage}
\hfill
\begin{minipage}[c]{0.49\linewidth}
\centering
\includegraphics[width=\textwidth]{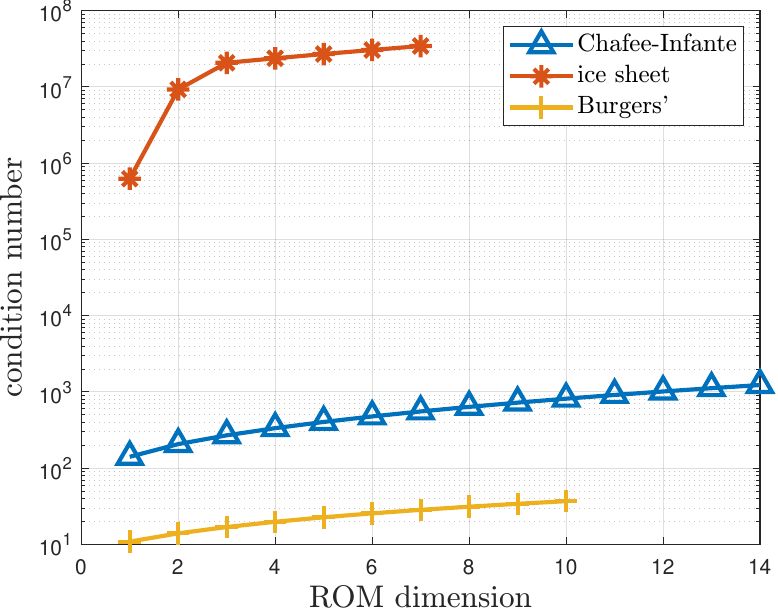}
\caption{\revfour{
Condition numbers of data matrices $\Pm$.}}
\label{fig:condition numbers}
\end{minipage}%
\end{figure}

\subsection{Chafee-Infante equation / Allen-Cahn equation}
\label{sec:chafee infante equation}

We consider the 
Chafee-Infante equation / Allen-Cahn equation
in a setup based on \cite{peherstorfer2020sampling,sawant2023physics},
\begin{align}
    \ppt \xs(\xis,t) = \frac{\partial^2}{\partial\xis^2} \xs(\xis,t) - \xs^3(\xis,t) + \xs(\xis,t),
    \label{eq:chafee-infante}
\end{align}
with spatial coordinate $\xis\in\Omega=(0,1)$ and time  
$t\in[0,0.1]$. 
We impose the boundary conditions
\begin{align}
    \xs(0,t) = \uus(t), \h 
    \ppx{\revfour{\xi}}
    \xs(1,t) = 0, \h\h t\in[0,0.1],
    \label{eq:chafee infante bc}
\end{align}
and zero initial condition $\xs(\xis,0) = 0$ for $\xis\in\Omega$. We discretize with central finite differences on an equidistant grid in space with mesh width $\Delta\xis = 2^{-7}$, and with explicit Euler in time with time step size $\Delta t=10^{-5}$, to obtain a time-discrete dynamical system of dimension $N=128$,
\begin{align}
    \xv_{k+1} = \xv_k + \Delta t\left(\Am_1 \xv_k + \Am_2\xv_k^2 + \Am_3 \xv_k^3 + \Bm \uus_k\right), \h k=0,\dots,K-1.
    \label{eq:chafee-infante discrete}
\end{align}
This system is polynomial of degree $l = 3$ with $\Am_1\in\mbR^{N\times N}$, $\Am_2\in\mbR^{N\times N_2}$ with $N_2 = \frac{N(N+1)} 2$, $\Am_3\in\mbR^{N\times N_3}$ with $N_3 = \frac{N(N+1)(N+2)} 6$, $\Bm\in\mbR^{N\times \Nu}$ with $\Nu = 1$, the input signal
$\uus_k=10(\sin(k \pi \Delta t) +1)$, $k=0,\dots,K-1$
and $K = 10^4$.  
The resulting 
$10^4 + 1$ 
snapshots are used to compute the POD bases
$\Vn\in\mbR^{N\times n}$ with $n=1,\dots,14$.

We apply Algorithm \ref{alg:exact opinf} with
\R{4: dt chafee}
$\Delta t = 
\revfour{3.2741}
\cdot 10^{-5}$
according to estimate \eqref{eq:dt estimate}
and compute the relative operator errors for ROM dimensions $n=1,\dots,14$.
Fig.\ \ref{fig:operator errors} shows that the operator error is on the order of machine precision and Fig. \ref{fig:condition numbers} shows that the data matrices $\Pm$ have moderate condition numbers.

A crucial advantage of our method is the minimal number of FOM time steps needed for exact reconstruction of the intrusive ROM. For ROM dimension $n=14$, we only need $\nf = n_1 + n_2 + n_3 + \Nu = 680$ FOM time steps.
This number is orders of magnitude smaller than in previous works with similar setups.
For example, in \cite{peherstorfer2020sampling}, 
$10^7$ FOM time steps with \reprojection\ have been used to reconstruct the operators exactly for ROM dimension $n=12$.
In \cite{sawant2023physics}, 
$10^5$ FOM time steps are used to infer ROM operators for ROM dimension $n=14$ from a regularized operator inference minimization problem that significantly differ from the intrusive counterparts.

Note that $\Am_2 = \mat 0$ in our discretization, so we could omit this term in the operator inference, so we would even need $n_2 = 105$ less FOM time steps.
We refrain from this simplification to keep the comparison with \cite{peherstorfer2020sampling} and \cite{sawant2023physics} fair.

\subsection{Ice sheet model}
\label{sec:ice sheet model}

We consider the shallow ice equations in a setup based on \cite{aretz2024enforcing},
\begin{align}
\revtwo{
    \ppt x(\xi,t) = 
    c_1
    \revfour{\ppx \xi\Bigg(}
    x^2(\xi,t) \ppx \xi x(\xi,t) 
    \revfour{\Bigg)}
    + 
    c_2
    \revfour{\ppx \xi\Bigg(}
    x^5(\xi,t) \left|\ppx \xi x(\xi,t)\right|^2 \ppx \xi x(\xi,t) 
    \revfour{\Bigg)}
    }
    \label{eq:ice sheet model}
\end{align}
with the spatial coordinate $\xi\in[0,1000]$, time $t\in[0,2]$ and the constants $c_1 = 
8.9
\cdot 10^{-13}$ and $c_2 = 
2.8
\cdot10^7$. We impose the initial condition
\begin{align}
    x(\xi,0) = 10^{-2} + 630\left(\frac \xi {2000} + 0.25\right)^4\left(\frac \xi {2000} - 0.75\right)^4,
\end{align}
and homogeneous Neumann boundary conditions. We discretize with finite differences in space with mesh width $\Delta \xi = \frac{1000}{2^9}$, and with implicit Euler in time with time step size $\Delta t = 10^{-3}$, to obtain the time-discrete dynamical system of dimension $N=512$,
\begin{align}
    \xv_{k+1} = \xv_k + \Delta t\left( \Am_3 \xv_k^3 + \Am_8 \xv_k^8 
    \right), \h k = 0,\dots,K-1.
\end{align}
This system is polynomial of degree $l=8$ with $\Am_3\in\mbR^{N\times N_3}$ with $N_3 = \frac{N(N+1)(N+2)} 6$ and $\Am_8\in\mbR^{N\times N_8}$ with $N_8 = {N+7 \choose 8}$ and $K=2000$. The resulting $2001$ snapshots from $t=0$ to $t=2$ are used to compute the POD bases $\Vn\in\mbR^{N\times n}$ with $n=1,\dots,7$.

We apply Algorithm \ref{alg:exact opinf} with
\R{4: dt ice}
$\Delta t = 
\revfour{1.4773}
\cdot 10^{13}$ according to estimate \eqref{eq:dt estimate}. This extreme value can be explained by the extreme operator norm $\|\tilde \Om\|_2 = 4.271\cdot10^{-14}$ for ROM dimension $n=1$. 
The relative operator errors for ROM dimensions $n=1,\dots,7$  are depicted in Fig.\ \ref{fig:operator errors}.
Like in the previous test case,
these operator errors are 
essentially on the order of machine precision. The fact that these errors are a few orders of magnitude larger 
than the errors in the previous test case is presumably related to the larger condition numbers as depicted in Fig.\ \ref{fig:condition numbers}.

For ROM dimension $n=7$, the proposed method performs $\nf=3087$ FOM time steps. This number is larger than the 
$2001$ 
FOM time steps with \reprojection\ reported in \cite{aretz2024enforcing}. However, the authors also report that this budget of snapshots is depleted before the inference of the ROM operators of dimension $n=6$ is completed. Hence, the\revtwo{ir} ROM operators of dimensions $n=6$ and $n=7$ are partially computed with regularized operator inference and snapshot data polluted with non-Markovian terms. Consequently, the\revtwo{ir} inferred operators of dimension $n=6$ and $n=7$ are not reconstructed exactly. The\revtwo{ir} operators of dimension $n\leq 5$ might be reconstructed exactly, but the used $2001$ snapshots are more than the $\nf = 365$ snapshots used by our method for dimension $n=5$.

Note that the cubic term could be neglected in practice due to the extremely small value of $c_1$. This simplification would reduce the number of required FOM time steps further.

\subsection{Burgers' equation}
\label{sec:burgers equation}

\begin{figure}
\centering
\includegraphics[width=.49\textwidth]{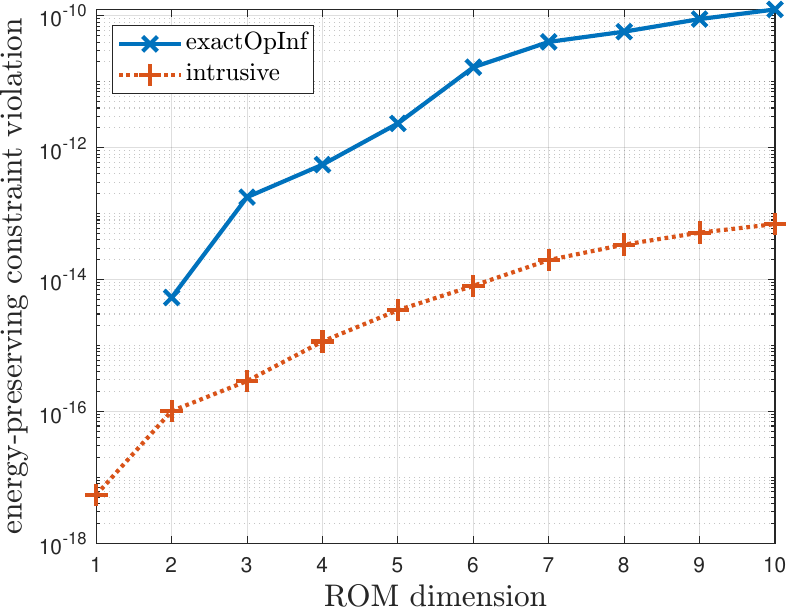}
\caption{\revfour{Violation of the energy-preserving constraint 
of the convection operator of Burgers' equation.}}
\label{fig:energy violation}
\end{figure}

\begin{figure}
\centering
\begin{minipage}[c]{0.49\linewidth}
\centering
\includegraphics[width=\textwidth]{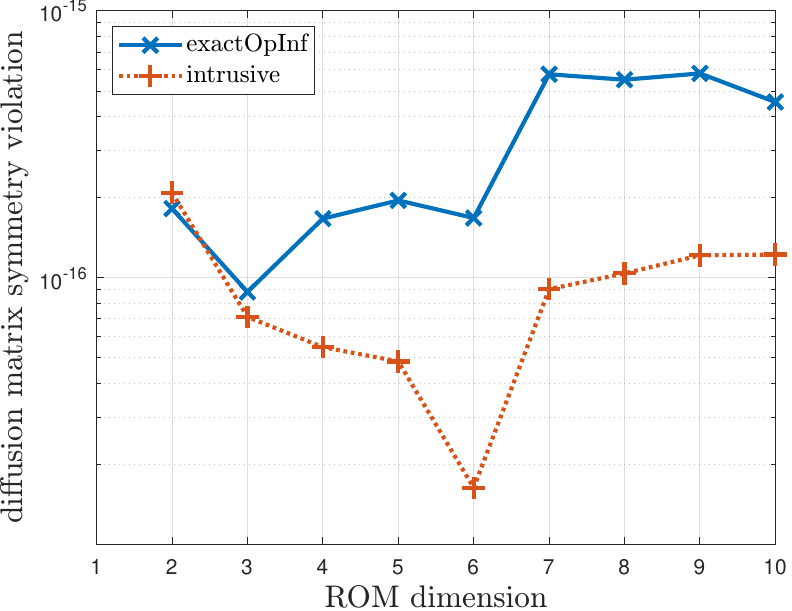}
\caption{\revfour{Symmetry violation of the diffusion matrix of Burgers' equation.}}
\label{fig:symmetry violation}
\end{minipage}
\hfill
\begin{minipage}[c]{0.49\linewidth}
\centering
\includegraphics[width=\textwidth]{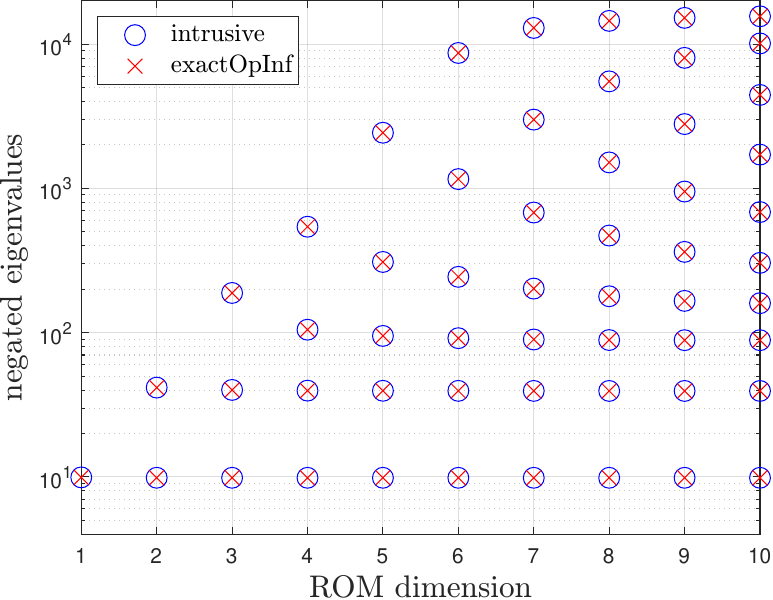}
\caption{\revfour{Eigenvalues of diffusion matrix of Burgers' equation 
multiplied with $-1$.}}
\label{fig:eigvals burgers periodic}
\end{minipage}%
\end{figure}

We consider the viscous Burgers' equation in a setup based on \cite{peherstorfer2020sampling,Koike_Qian_2024},
\begin{align}
    \ppt \xs(\xis,t) + \xs(\xis,t)\ppx{\xis}\xs(\xis,t) - 
    \frac{\partial^2}{\partial \xis^2}\xs(\xis,t) = 0,
    \label{eq:burgers equation}
\end{align}
\revtwo{
with 
spatial} coordinate $\xis\in\Omega=(-1,1)$, time $t\in[0,1]$.
We impose 
periodic boundary conditions and the initial condition
\R{4: burgers IC}
$
    x(\xi,0) = -\sin\left(
    \revfour{\pi}\xi\right)
    $
for all $\xi\in\Omega$.

We discretize with finite differences on an equidistant grid in space with mesh width $\Delta\xis=2^{-6}$, and with explicit Euler in time with step size $\Delta t= 10^{-4}$, to obtain the time-discrete dynamical system of dimension $N=128$,
\begin{align}
    \xv_{k+1} = \xv_k + \Delta t\left(\Am_1 \xv_k + \Am_2 \xv_k^2 
    \right), 
    \h k = 0,\dots,K-1.
\end{align}
This system is polynomial of degree $l=2$ with $\Am_1\in\mbR^{N\times N}$, $\Am_2\in\mbR^{N\times N_2}$ with $N_2 = \frac{N(N+1)}{2}$  
and $K=10^4$. 
Concretely, we use a spatial discretization as described in \cite{van2024energy} such that the diffusion matrix $\Am_1$ is symmetric and negative semi-definite, and the convection operator $\Am_2$ satisfies the skew-symmetry property $\xv^T \Am_2 \xv^2 = 0$ for all $\xv\in\mbR^N$. 
The resulting $10^4 + 1$ snapshots are used to compute the POD bases $\Vn\in\mbR^{N\times n}$ with $n=1,\dots,10$.

We perform Algorithm \ref{alg:exact opinf} 
with 
\R{4: dt burgers}
$\Delta t = 
\revfour{
1.0134\cdot 10^{-1}}
$ according to estimate \eqref{eq:dt estimate}
and compute the relative operator errors for ROM dimensions $n=1,\dots,10$.
Fig.\ \ref{fig:operator errors} shows operator errors on the order of machine precison and Fig.\ \ref{fig:condition numbers} very small condition numbers.
In addition, we 
observe that both, the intrusive and
the inferred operators preserve the 
skew-symmetry of the FOM convection operator and the symmetry and negative semi-definiteness of the FOM diffusion operator.
Fig.\ \ref{fig:energy violation} shows the energy-preserving constraint violation
$\sum_{1\leq i,j,k\leq n} \big|\hat h_{ijk} + \hat h_{jik} + \hat h_{kji}\big|$ as defined in \cite{Koike_Qian_2024} for the inferred operator with $\hat h_{ijk} = \onehot_i \hat \Am_2 \I_N^{(2)}(\onehot_j\otimes \onehot_k)$ and for the intrusive operator with $\hat h_{ijk} = \onehot_i \tilde \Am_2 \I_N^{(2)}(\onehot_j\otimes \onehot_k)$. As can be seen, this violation is essentially on the order of machine precision for both, inferred and intrusive operators.
Fig.\ \ref{fig:symmetry violation} shows the symmetry violation of the inferred and the intrusive diffusion matrix, $\frac{\left\|\hat \Am_1 - \hat \Am_1^T \right\|}{\left\|\hat \Am_1 \right\|}$ and $\frac{\left\|\tilde \Am_1 - \tilde \Am_1^T \right\|}{\left\|\tilde \Am_1 \right\|}$, respectively. Also this violation is on the order of machine precision. Consequently, the diffusion matrices are symmetric, so their eigenvalues are real. These eigenvalues - multiplied with $-1$ - are depicted in Fig.\ \ref{fig:eigvals burgers periodic}. As can be seen the eigenvalues of the inferred diffusion matrices coincide with the eigenvalues of the corresponding intrusive matrices. Hence, both ROMs in particular preserve the negative semi-definiteness of the FOM diffusion operator.


\section{Conclusion}
\label{sec:conclusion}

This work 
introduces
a method to generate 
rank-ensuring
snapshot data for operator inference that 
guarantees
the exact reconstruction of the corresponding intrusive reduced order model. 
This method performs the minimal number of FOM time steps and yields a square data matrix, so the operator inference least squares problem can be simplified to a linear system that has better numerical stability.
Numerical experiments 
confirm
that the novel snapshot data generation method 
leads to exact reconstruction of the intrusive ROM operators. Consequently, also 
properties
such as skew-symmetry, symmetry and negative semi-definiteness of the intrusive operators are preserved by the \nonintrusive\ operators without the need to enforce this structure preservation explicitly.

To improve the condition number of the matrix corresponding to the operator inference linear system, future work could investigate
the exploitation of nestedness to solve subproblems sequentially \cite{aretz2024enforcing} or different choices of rank-ensuring ROM states, for example based on Chebychev points
\cite{trefethen2019approximation}.
\R{1:future work}\revone{Future work should also investigate how to 
enable that the proposed rank-ensuring states can be specified as initial conditions in solvers that allow only a certain range of state values or require specific boundary conditions.}

\newcommand{\Wn}{\mat W_n}

In this work, we have not considered the parameter dependency of ROM operators.
For that purpose,
our work could be combined with \cite{mcquarrie2023nonintrusive} to generate rank-sufficient operator inference snapshot data for parametric problems. We expect the full rank guarantee to be transferable to this setting.

\appendix

\input{sections/existence_and_uniqueness}

\revtwo{
\section{Algorithm for generating ROM state set}
\label{appendix: algo gen rom states}
We present two different approaches how to implement the procedure to generate the rank-ensuring ROM state sets $\setXi$ as defined in \eqref{eq:def set X i} in line \ref{alg line:gen rom state set} of Algorithm \ref{alg:exact opinf}.
As shown in Algorithm \ref{alg:gen rom state set 3} for polynomial degree $i=3$, the procedure essentially comprises 3 nested \textbf{for}-loops. The first approach shown in Algorithm \ref{alg:gen rom state set recursive} constructs these nested \textbf{for}-loops for arbitrary $i$ via recursion. This recursion makes the algorithm rather non-intuitive. While this non-intuitivity might be unavoidable to handle arbitrary polynomial degrees $i$, a more intuitive approach is possible in practice. In practice, we know a priori which polynomial degrees $i$ are involved in the considered problems. For example, in the numerical experiments in Section \ref{sec:num exps}, the only polynomial degrees involved are $\{1,2,3,8\}$. Consequently, we can use the hard-coded discrimination of cases shown in Algorithm \ref{alg:gen rom state set - hard coded} as alternative to Algorithm \ref{alg:gen rom state set recursive}, with the procedures for $i=1,2,8$ implemented analogous to Algorithm \ref{alg:gen rom state set 3}. 
}

\begin{algorithm}
\caption{
\revtwo{
Generate rank-ensuring ROM state set $\setXx 3$
}
}
\label{alg:gen rom state set 3}
\revtwo{
\begin{algorithmic}[1]
\Procedure{GenerateROMstateSet3}
{
ROM dimension $n$
    }
    \State Initialize $\setXx 3 \gets \emptyset$
    \For{$j_1 = 1,\dots,n$}
    \For{$j_2 = j_1,\dots,n$}
    \For{$j_3 = j_2,\dots,n$}
        \State $\setXx 3 \gets \setXx 3 \cup \onehot_{j_1} + \onehot_{j_2} + \onehot_{j_3}$
    \EndFor
    \EndFor
    \EndFor
    \State
\Return $\setXx 3$
\EndProcedure
\end{algorithmic}
}
\end{algorithm}

\newcommand{\fordepth}{{i^*}}

\begin{algorithm}
\caption{
\revtwo{
Generate rank-ensuring ROM state set $\setXi$ - recursive algorithm
}
}
\label{alg:gen rom state set recursive}
\revtwo{
\begin{algorithmic}[1]
\Procedure{GenerateROMstateSet}
{
\newline\phantom{---} ROM dimension $n$,
\newline\phantom{---} polynomial degree $i\in\mbN_0$
    \newline}
    \State \Return \Call{recursiveROMstateSet}{
    ROM dimension $n$, polynomial degree $i$, 
    \newline \phantom{-------------------------------------------------------} current for loop depth $\fordepth=0$, 
    \newline \phantom{-------------------------------------------------------} current index sequence $\vt j = []$, 
    \newline \phantom{-------------------------------------------------------} current ROM state set $\setXi=\emptyset$
    }
\EndProcedure
\newline
\Procedure{recursiveROMstateSet}
{
    ROM dimension $n$, polynomial degree $i\in\mbN_0$, 
    \newline \phantom{----------------------------------} current for loop depth $\fordepth$, 
    \newline \phantom{----------------------------------} current index sequence $\vt j = [j_1, \dots,j_\fordepth]$, \Comment{always has length $\fordepth$} 
    \newline \phantom{----------------------------------} current ROM state set $\setXi$
    }
    \If{$\fordepth = i$}
    \State Initialize $\cx \gets \vt 0 \h \in\mbR^n$
    \For{k=1,\dots,i}
        \State $\cx \gets \cx + \onehot_{j_k}$
    \EndFor
    \State \Return $\setXi \cup \cx$
    \Else
    \LineComment{Perform next inner for-loop}
    \For{$j_{\fordepth+1} = j_\fordepth,\dots,n$} \Comment{with $j_0 = 1$}
    \State \Return \Call{recursiveROMstateSet}{
    ROM dimension $n$, polynomial degree $i$, 
    \newline \phantom{-------------------------------------------------------} current for loop depth $\fordepth+1$, 
    \newline \phantom{-------------------------------------------------------} current index sequence $[j_1,\dots,j_{\fordepth+1}]$, 
    \newline \phantom{-------------------------------------------------------} current ROM state set $\setXi$   
    }
    \EndFor
    \EndIf
\EndProcedure
\end{algorithmic}
}
\end{algorithm}


\begin{algorithm}
\caption{
\revtwo{
Generate rank-ensuring ROM state set $\setXi$ - hard-coded
}
}
\label{alg:gen rom state set - hard coded}
\revtwo{
\begin{algorithmic}[1]
\Procedure{GenerateROMstateSet}
{
\newline\phantom{---} ROM dimension $n$,
\newline\phantom{---} polynomial degree $i\in\mbN_0$
    \newline}
    \If{$i=1$} \Return \Call{GenerateROMstateSet1}{ROM dimension $n$}
    \ElsIf{$i=2$} \Return \Call{GenerateROMstateSet2}{ROM dimension $n$}
    \ElsIf{$i=3$} \Return \Call{GenerateROMstateSet3}{ROM dimension $n$}
    \ElsIf{$i=8$} \Return \Call{GenerateROMstateSet8}{ROM dimension $n$}
    \EndIf
\EndProcedure
\end{algorithmic}
}
\end{algorithm}

\begin{figure}
\centering
\begin{minipage}[c]{0.49\linewidth}
\centering
\includegraphics[width=\textwidth]{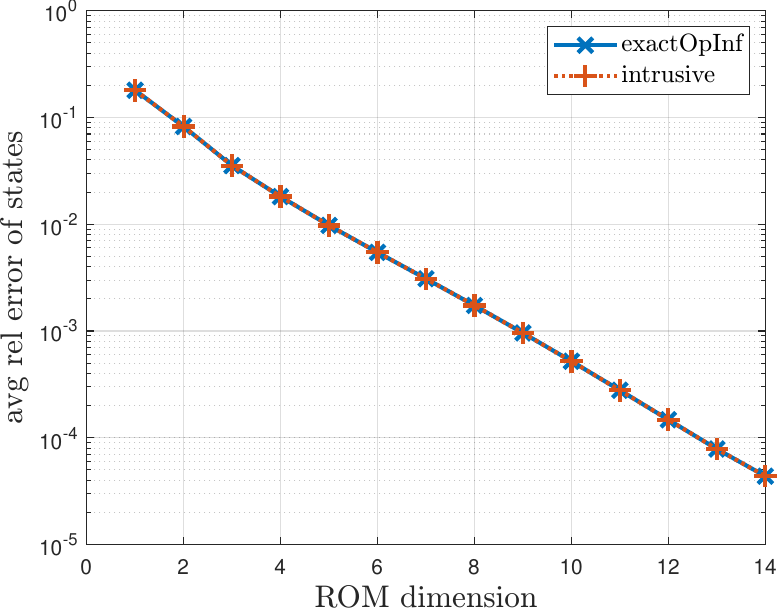}
\caption{
\revtwo{
Relative average ROM state error \eqref{eq: rel avg rom state error} for Chafee-Infante equation.}
}
\label{fig:rom state error chafee}
\end{minipage}
\hfill
\begin{minipage}[c]{0.49\linewidth}
\centering
\includegraphics[width=\textwidth]{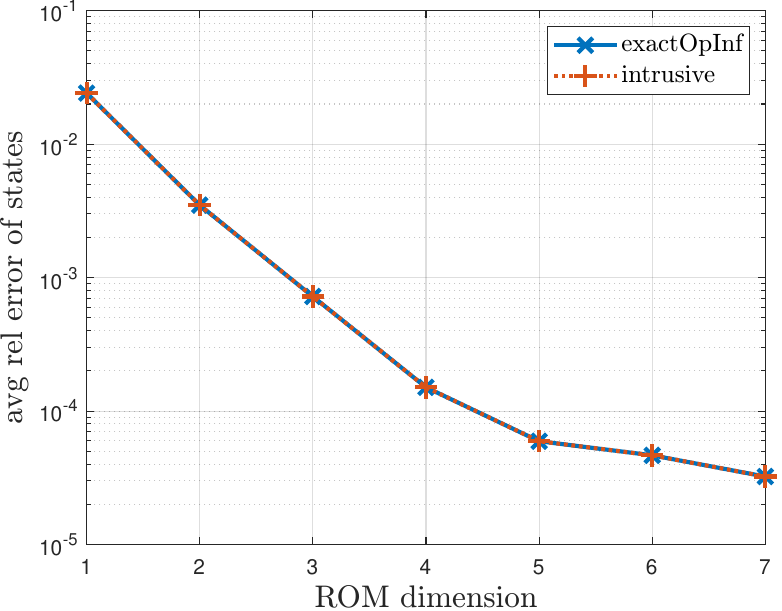}
\caption{\revtwo{
Relative average ROM state error \eqref{eq: rel avg rom state error} for the ice sheet model.}}
\label{fig:rom state error ice}
\end{minipage}%
\end{figure}

\begin{figure}
\centering
\includegraphics[width=.49\textwidth]{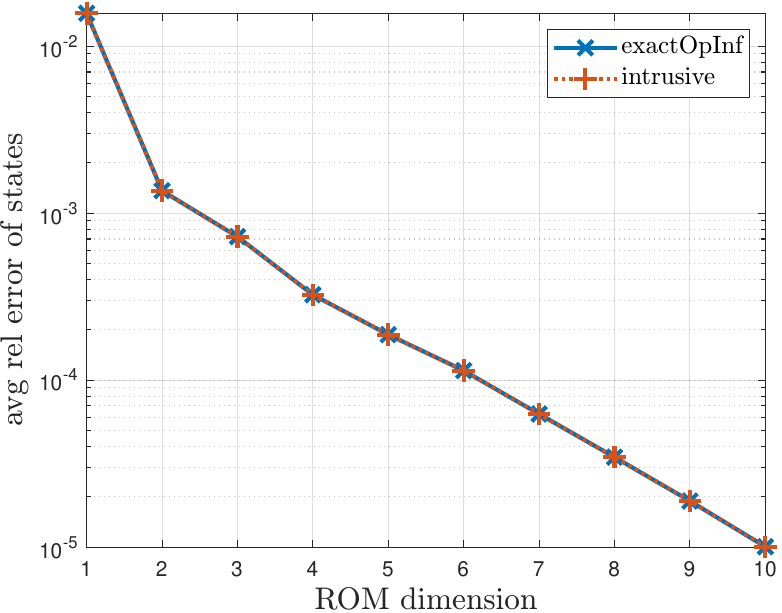}
\caption{\revtwo{Relative average ROM state error \eqref{eq: rel avg rom state error} for Burgers' equation.}}
\label{fig:rom state error burgers}
\end{figure}

\revtwo{
\section{ROM state errors with respect to the FOM}
\label{appendix:rom state errors}
Figures \ref{fig:rom state error chafee}, \ref{fig:rom state error ice} and \ref{fig:rom state error burgers} show the relative average ROM state errors,
\begin{align}
    \frac{\left\|\mat X - \Vn \mat{\tilde X}\right\|_F}{\left\|\mat X\right\|_F},
    \label{eq: rel avg rom state error}
\end{align}
for Chafee-Infante equation in Section \ref{sec:chafee infante equation}, ice sheet model in Section \ref{sec:ice sheet model} and Burgers' equation in Section \ref{sec:burgers equation}, respectively.
The matrix $\mat X = \left[\x_0 \;\; \dots \;\; \x_{K_\text{POD}-1}\right]$ comprises the FOM snapshots used to construct the POD basis, and $\mat{\tilde X} = \left[\tx_0 \;\; \dots \;\; \tx_{K_\text{POD}-1}\right]$ the corresponding ROM snapshots generated by integrating the ROM ODE \eqref{eq:rom efficient} with the intrusive ROM operators and the inferred ROM operators, respectively. 
As can be seen, the ROM errors decrease monotonously for increasing ROM dimensions in all three test cases for both, intrusive and inferred ROMs. Additionally, we cannot see any difference between the errors of corresponding intrusive and inferred ROMs. In fact, the differences between these errors are all smaller than $10^{-10}$.
This alignment is due to the exact reconstruction of the intrusive ROM operators
shown in Fig. \ref{fig:operator errors}.
} 

\bibliographystyle{siamplain}
\bibliography{library/references}

\end{document}

%% file: ex_shared.tex
\usepackage{lipsum}
\usepackage{amsfonts}
\usepackage{graphicx}
\usepackage{epstopdf}
\ifpdf
  \DeclareGraphicsExtensions{.eps,.pdf,.png,.jpg}
\else
  \DeclareGraphicsExtensions{.eps}
\fi

\newsiamremark{remark}{Remark}
\newsiamremark{hypothesis}{Hypothesis}
\crefname{hypothesis}{Hypothesis}{Hypotheses}
\newsiamthm{claim}{Claim}
\newsiamremark{fact}{Fact}
\crefname{fact}{Fact}{Facts}

\headers{Exact operator inference with minimal 
data}
{H. Rosenberger, B. Sanderse, and G. Stabile}

\title{Exact operator inference with minimal 
data\thanks{Submitted to the journal's Methods and Algorithms for Scientific Computing section  June 27, 2025.}
\funding{This publication is part of the project “Discretize first, reduce next” (with project number VI.Vidi.193.105 of the research
programme NWO Talent Programme Vidi which is (partly) financed by the Dutch Research Council (NWO).}}

\author{Henrik Rosenberger\footnotemark[2] \footnotemark[3]\thanks{Centrum Wiskunde \& Informatica, Amsterdam, The Netherlands 
  (\email{henrik.rosenberger@cwi.nl} \email{sanderse@cwi.nl}).}
\and Benjamin Sanderse\footnotemark[2] \footnotemark[3]\thanks{Centre for Analysis, Scientific Computing and Applications, Eindhoven University of Technology, Eindhoven, The Netherlands.}
\and Giovanni Stabile\thanks{The Biorobotics Institute, Sant’Anna School of Advanced Studies, Pisa, Italy
  (\email{giovanni.stabile@santannapisa.it}).}}

\usepackage{amsopn}


%% file: preamble/packages.tex
\usepackage[utf8]{inputenc}
\usepackage[T1]{fontenc}
\usepackage{ae}
\usepackage[english]{babel}
\usepackage{amsmath,amssymb,amstext,amsfonts}
\usepackage{bm}
\usepackage{graphicx}
\usepackage{pgf}
\usepackage{pgfplots}

\usepackage{tikz}
\usetikzlibrary{arrows}
\usetikzlibrary{shapes}
\usetikzlibrary{decorations.text}
\usepackage{xcolor}
\usepackage{color}
\usepackage{tikzscale}
\usepackage{pgfplots}

\usetikzlibrary{shapes.geometric}
\usetikzlibrary{positioning}
\usetikzlibrary{positioning,shapes.multipart, fit,backgrounds,calc}

\usepackage{nicefrac}
\usepackage{array,tabularx}
\usepackage{multirow}
\newcolumntype{L}[1]{>{\raggedright\arraybackslash}p{#1}}
\newcolumntype{C}[1]{>{\centering\arraybackslash}p{#1}}
\newcolumntype{R}[1]{>{\raggedleft\arraybackslash}p{#1}}
\usepackage{booktabs}
\usepackage{mathtools}

\usepackage{geometry}
\usepackage{pdfpages}

\usepackage[colorinlistoftodos,prependcaption,textsize=tiny]{todonotes}
\usepackage{caption}
\usepackage{subcaption}
\usepackage{framed}

\usepackage{siunitx} 
\usepackage[normalem]{ulem} 

\theoremstyle{break}
\newtheorem{example}{Example}

\usepackage{algorithm}
\usepackage{algpseudocode}

\usepackage[nocompress]{cite} 

%% file: preamble/macros.tex
\newcommand{\vt}[1]{{\bm{#1}}}
\newcommand{\mat}[1]{{\mathbf{#1}}}

\newcommand{\ppt}{\frac{\partial}{\partial t}}

\newcommand{\ppx}[1]{\frac{\partial}{\partial #1}}

\newcommand{\mbR}{\mathbb{R}}
\newcommand{\mbN}{\mathbb{N}}

\newcommand{\h}{\hspace{1 em}}

%% file: preamble/rank_suff_macros.tex
\newcommand{\x}{\vt x}
\newcommand{\xv}{\vt x}
\newcommand{\uv}{\vt u}
\newcommand{\f}{\vt f}
\newcommand{\Am}{\mat A}
\newcommand{\B}{\mat B}
\newcommand{\Bm}{\mat B}

\newcommand{\I}{\mat I}
\newcommand{\Jm}{\mat J}

\newcommand{\tx}{\tilde \x}
\newcommand{\Vn}{\mat V_n}
\newcommand{\vv}{\vt v}

\newcommand{\tk}{t_k} 
\newcommand{\bx}{\breve \x}

\newcommand{\Om}{\mat O}
\newcommand{\pv}{\vt p}

\newcommand{\mdots}{\;\;\dots\;\;}

\newcommand{\Pm}{\mat P}
\newcommand{\bX}{\breve{\mat X}}

\newcommand{\onehot}{\vt{e}}

\newcommand{\nf}{{n_f}}
\newcommand{\np}{{n_p}}

\newcommand{\ind}{s}

\newcommand{\inds}{\mathcal I}

\newcommand{\muv}{\vt \mu}
\newcommand{\muj}{\muv}

\newcommand{\sumi}{\sum_{i\in\inds}}

\newcommand{\reprojection}{re-projection}

\newcommand{\degreeset}{degree set}
\newcommand{\nonintrusive}{nonintrusive}

\newcommand{\Nu}{{N_u}}

\newcommand{\setXi}{\mathcal X^i}
\newcommand{\Um}{\mat U}
\newcommand{\sx}{\tx}

\newcommand{\poly}{p}

\newcommand{\xcoff}{c}

\newcommand{\ls}{l}
\newcommand{\Nl}{\mbN_\ls}
\newcommand{\Is}{\mathcal I}
\newcommand{\PI}{\mathcal{P}^1_\Is}

\newcommand{\cs}{c}
\newcommand{\cv}{\vt c}
\newcommand{\fv}{\vt f}
\newcommand{\ival}{f}
\newcommand{\nI}{{n_\Is}}

\newcommand{\polyii}{\bar\poly}
\newcommand{\intmat}{\mat A}

\newcommand{\setXx}[1]{\mathcal X^{#1}}

\newcommand{\alphav}{\vt \alpha}
\newcommand{\Xv}{\vt y}
\newcommand{\qv}{\vt q}

\newcommand{\n}{{m}}

\newcommand{\polyi}{\hat\poly}
\newcommand{\yv}{\hat \xv}
\newcommand{\betav}{{\hat \alphav}}
\newcommand{\hatqv}{\hat \qv}

\newcommand{\map}{\vt m}

\newcommand{\pminus}{p^-}
\newcommand{\pplus}{p}
\newcommand{\indsminus}{{\inds^-}}
\newcommand{\indsplus}{\inds}
\newcommand{\inew}{{i^*}}
\newcommand{\indsmi}{{\inds^{-\inew}}}

%% file: preamble/main_only_macros.tex
\algnewcommand{\LineComment}[1]{\vspace{.0625in}\State \textcolor{gray}{\texttt{\# #1}}}
\algrenewcomment[1]{\hfill\textcolor{gray}{\texttt{\# #1}}}

\newtheorem{mytheo}[theorem]{Theorem}

%% file: revision/report_commands.tex
\input{revision/rev_commands}

\newcommand*\linenomathpatchAMS[1]{%
  \expandafter\pretocmd\csname #1\endcsname {\linenomathAMS}{}{}%
  \expandafter\pretocmd\csname #1*\endcsname{\linenomathAMS}{}{}%
  \expandafter\apptocmd\csname end#1\endcsname {\endlinenomath}{}{}%
  \expandafter\apptocmd\csname end#1*\endcsname{\endlinenomath}{}{}%
}
\expandafter\ifx\linenomath\linenomathWithnumbers
  \let\linenomathAMS\linenomathWithnumbers
  \patchcmd\linenomathAMS{\advance\postdisplaypenalty\linenopenalty}{}{}{}
\else
  \let\linenomathAMS\linenomathNonumbers
\fi
\linenomathpatchAMS{gather}
\linenomathpatchAMS{multline}
\linenomathpatchAMS{align}
\linenomathpatchAMS{alignat}
\linenomathpatchAMS{flalign}


%% file: revision/rev_commands.tex
\newcommand{\R}[1]{}

\newcommand{\revone}[1]{{#1}}
\newcommand{\revtwo}[1]{{#1}}
\newcommand{\revthree}[1]{{#1}}
\newcommand{\revfour}[1]{{#1}}

%% file: sections/introduction.tex
\section{Introduction}

{Computational efficiency is of paramount importance in complex engineering tasks such as optimization, data assimilation, uncertainty quantification and the development of digital twins.} All these tasks require multiple simulations with similar settings. While first principle-derived models have primarily been developed to run single simulations, data-driven models are being developed to exploit the similarity among such multiple simulations. 

{These data-driven surrogate models vary in how much they preserve structure of the corresponding first principle-derived model \cite{kramer2024learning}.} 
Statistical data-fit models such as Gaussian process models \cite{williams2006gaussian}
generally do not preserve such structure and
solely approximate the mapping between simulation inputs and outputs.
Reduced-order models (ROMs), on the other hand, are closely based on first principle-derived, full-order models (FOMs) and perform simulations on a lower-dimensional subspace, for example, by projecting the first principle-derived model onto this subspace \cite{antoulas2005approximation,rozza2008reduced,hesthaven2016certified}. Advantages of such projection-based ROMs are the availability of rigorous theory that provides a posteriori error estimates \cite{prud2002reliable,veroy2005certified,hinze2005proper,haasdonk2011efficient,urban2012new}, and the smaller demand for training data compared with 
surrogate models
that preserve less structure. 
A disadvantage of projection-based ROMs is that they are intrusive in the sense that they require access to the numerical operators of the 
FOM 
to perform the projection onto the lower-dimensional subspace. This access is often not possible in software used in industry, 
and hence impedes
the use of intrusive reduced-order models in practice.

{To address this practical problem, methods have been developed to approximate intrusive reduced models in a \nonintrusive\ way.} These \nonintrusive\ methods do not require access to the numerical operators of the 
FOM,
but only to the generated output data. 
For linear models, dynamic mode decomposition finds the linear operator that fits the data best \cite{kutz2016dynamic}.
For nonlinear models, operator inference generalizes this approach to 
polynomial operators
\cite{PEHERSTORFER2016data}. For non-polynomial models, a variable transformation can be used to obtain a polynomial model \cite{qian2020lift}, and so operator inference can still be applied.

However, two major difficulties currently impede operator inference in accurately reconstructing the intrusive ROM operators.
First, the FOM snapshot data contains non-Markovian information that is inconsistent with the Markovian intrusive ROM.
This non-Markovian information \revone{perturbs} the operator inference minimization problem such that the intrusive ROM generally does not minimize this perturbed problem.
To eliminate this  
perturbation,
Peherstorfer has proposed to generate additional FOM snapshot data whereby the FOM state is projected onto the ROM basis after each time step \cite{peherstorfer2020sampling}, known as the \reprojection\ approach. The resulting FOM snapshot data does not contain non-Markovian information and is equivalent to snapshot data 
generated by the intrusive ROM. 
Consequently, 
the minimization problem is unperturbed and the intrusive ROM does minimize it. 
However, operator inference is only guaranteed to yield this solution if this 
minimizer 
is unique. Concretely, the operator inference minimization problem is a least squares problem and has a unique solution only if the least squares matrix has full rank.

This full rank requirement is the second major difficulty in accurately reconstructing the intrusive ROM operators.
To the authors' knowledge, no method has been proposed yet to ensure this full rank requirement a priori. In practice, two different approaches have been suggested to address this issue. The first approach is to generate large numbers of snapshots and test a posteriori whether the resulting least squares problem has full rank \cite{PEHERSTORFER2016data,peherstorfer2020sampling}. However, this approach 
does not guarantee full rank, even when
the number of generated snapshots exceeds the required rank by orders of magnitude, and hence causes large computational costs. 
The second approach is to circumvent the need for full rank by adding a regularization term to the operator inference least squares problem \cite{swischuk2020learning,sawant2023physics}. However, this modification inevitably causes the resulting 
inferred ROM operators
to deviate from the intrusive ROM operators.

In this work, we propose a novel snapshot data generation method that guarantees the full rank of the operator inference least squares problem. 
The number of generated snapshots is equal to this rank and hence minimal. Like Peherstorfer's approach \cite{peherstorfer2020sampling}, our method avoids non-Markovian information and hence 
reconstructs the intrusive ROM operators 
exactly, i.e., up to machine precision.
Our method has two key ingredients.
First, we leverage the knowledge about the specific polynomial structure of the intrusive ROM to identify a set of 
ROM states
and input signals that define snapshots that induce a full rank least squares problem.
Second, we generate snapshots corresponding to these identified ROM states and input signals by 
simulating multiple trajectories, initialized with the identified ROM states as initial conditions and the identified input signals, for only a single 
\R{1:explicit euler intro 1}\revone{explicit Euler} time step.
In contrast to the usual approach that generates snapshot data along long trajectories, this single-step approach gives us 
control over which snapshot data we obtain. Additionally, this single-step approach avoids any non-Markovian terms, so we do not need \reprojection \cite{peherstorfer2020sampling}\R{1:explicit euler intro 2}\revone{, and the explicit Euler scheme prevents time discretization inconsistencies}.
Since the number of snapshots is minimal, the least squares problem simplifies to a linear system that is numerically more stable.
Since the method exactly reconstructs intrusive ROM operators, it preserves structures such as symmetry and skew-symmetry of these intrusive operators \cite{sanderse2020non} without the need to enforce this preservation explicitly \cite{Koike_Qian_2024,gruber2023canonical,goyal2023guaranteed}.
\R{1:no reproj}

This article is organized as follows.
In Section \ref{sec:preliminaries}, we introduce dynamical systems with polynomial terms, intrusive projection-based model reduction, operator inference and the limitation of insufficient snapshot data. In Section \ref{sec:rank suff snapshot data generation}, we propose our novel method for rank-sufficient snapshot data generation. In Section \ref{sec:num exps}, we demonstrate the novel method in numerical experiments.  
In Section \ref{sec:conclusion}, we draw conclusions and give an outlook on future work.

%% file: sections/preliminaries.tex
\section{Preliminaries}
\label{sec:preliminaries}

Our notation is based on \cite{peherstorfer2020sampling}, but we start at the time-continuous level.

\subsection{Dynamical systems with polynomial  
terms}
\label{sec:prelimi fom}
As full-order model (FOM), we
consider a time-continuous dynamical system 
\begin{align}
    \dot \x(t;\muv) = \f(\x(t;\muv),\uv(t;\muv);\muv), \h \x(0;\muv) = \x_0(\muv),
    \label{eq:peher fom time conti}
\end{align}
with state $\x(t;\muv)\in\mbR^N$, input $\uv(t;\muv)\in\mbR^\Nu$, time $t\in[0,T]$,
and parameter $\muv\in\mathcal{D} \subset \mbR^d$.
We assume the nonlinear function $\f:\mbR^N\times\mbR^\Nu\times\mathcal D \rightarrow \mbR^N$ to be polynomial of order $l\in\mbN_0$ in the state $\x$ and linear in the input signal $\uv$, so there exist $\Am_i(\muv)\in\mbR^{N\times N_i}$ with $N_i = \begin{pmatrix}
    N+i-1 \\ i
\end{pmatrix}$ for 
$i\in\inds$,
and $\B(\muv)\in\mbR^{N\times \Nu}$ for $\muv\in\mathcal D$ such that
\begin{align}
    \f(\x,\uv;\muv) = \sum_{i\in\inds} \Am_i(\muv)\x^i + \B(\muv)\uv,
    \label{eq:fom dynamics}
\end{align}
with the \degreeset\ $\inds\subset\{0,1,\dots,l\}$.
The vector $\x^i \in\mbR^{N_i}$ comprises all unique products of $i$ entries of $\x$ ($\x^0$ is scalar 1). We can construct this vector from the $i$-fold Kronecker product $\x\otimes \dots \otimes \x=: \x^{i_\otimes}$ by removing all duplicate entries due to commutativity of the multiplication \cite{PEHERSTORFER2016data}.
Formally, we can write 
$
    \x^i = \I_N^{(i)} 
    \x^{i_\otimes}
    $
with $\I_N^{(i)}\in\{0,1\}^{N_i\times N^i}$ an identity matrix with certain rows removed, and conversely
$
    \x^{i_\otimes} = \Jm_N^{(i)} \x^i
    $
with $\Jm_N^{(i)}\in\{0,1\}^{N^i\times N_i}$ an identity matrix with certain rows copied and permuted.
Here $N^i$ denotes as usual the $i$-th power of $N$.

\subsection{Intrusive projection-based model reduction}
\label{sec:prelimi intrusive rom}

To obtain a projection-based reduced order model (ROM), a common intrusive approach is to 
project the FOM \eqref{eq:peher fom time conti} onto $\Vn$ and
replace the FOM state $\x(t;\muv)$ by the lower-dimensional approximation $\Vn \tx(t;\muv)$,
\begin{align}
    \dot \tx(t;\muv) = \Vn^T \f(\Vn \tx(t;\muv),\uv(t;\muv);\muv), \h \tx(0;\muv) = \Vn^T \x_0(\muv). 
    \label{eq:rom naive}
\end{align}
Here, $\Vn = [\vv_1, \dots, \vv_n]\in\mbR^{N\times n}$ is an orthonormal basis of an $n$-dimensional subspace of $\mbR^N$ with $n\ll N$, and $\tx(t;\muv)\in\mbR^n$ the ROM state.
The orthonormal basis $\Vn$ can, for example, be obtained from proper orthogonal decomposition (POD)
\cite{hinze2005proper} 
of snapshots $\x_k(\muv)\in\mbR^N, \, k = 0,\dots,K_\text{POD}-1$
of a time discretization of the FOM \eqref{eq:peher fom time conti}. Other methods to construct reduced bases and details on how to sample parameters and the ROM dimension $n$
are discussed in \cite{benner2015survey}.

To allow efficient simulations of the ROM \eqref{eq:rom naive}, the cost of evaluating the right-hand side must be independent of the FOM dimension $N$. For this purpose, the intrusive approach is to precompute reduced operators,
    \begin{align} 
        \tilde \B(\muj) := \Vn^T \B(\muj) \h \in\mbR^{n\times \Nu}, 
        \label{eq:intrusive ops1}
    \end{align}
    and
        \begin{align}
        \tilde \Am_i(\muj) := \Vn^T \Am_i(\muj)\I_N^{(i)}
        \Vn^{i_\otimes}
        \Jm_N^{(i)} \; \in\mbR^{n\times n_i}, \h 
        \; \text{ with } \;
        n_i = \begin{pmatrix}
            n + i - 1 \\ i
        \end{pmatrix}, \; \text{ for all } \; i\in\inds,
        \label{eq:intrusive ops2}
    \end{align}
whose dimensions only depend on the ROM dimension $n$ (and the unreduced input signal dimension $\Nu$) \cite{peherstorfer2020sampling}. 
Note that these precomputations must in general be performed for each choice of the parameter $\muj\in\mathcal{D}$ separately.
The resulting efficient ROM for fixed $\muj$ is given by
    \begin{align}
        \dot \tx(t;\muj) 
        &= \tilde \f(\tx(t;\muj),\uv(t;\muj);\muj) 
        =
        \sumi
        \tilde \Am_i(\muj)\tx^i(t;\muj) + \tilde \B(\muj)\uv(t;\muj) 
        =
        \tilde \Om(\muj) \pv(\tx(t;\muj),\uv(t;\muj)),
        \label{eq:rom efficient} 
    \end{align}
with the aggregated ROM operator 
$\tilde \Om(\muv) := \left[\left[\tilde \Am_i(\muv)\right]_{i\in\inds} \;\; \tilde \Bm(\muv) \right] \;\in\mbR^{n \times \nf}$, 
$\nf = \np + \Nu$, $\np = \sumi n_i$ and the 
aggregated ROM state (and input) polynomial
\begin{align}
    \pv: \mbR^n \times \mbR^\Nu \rightarrow \mbR^\nf, \h 
    \pv(\tx,\uv) = \begin{bmatrix}
        \left[\tx^i\right]_{i\in\inds} \\ \uv
    \end{bmatrix}.
    \label{eq:def p vec}
\end{align}

\subsection{Operator inference}
\label{sec:standard opinf}

As described in the previous subsection, intrusive model reduction requires access to the FOM operators 
$\Am_i(\muv)$, $\Bm(\muv)$
to precompute the ROM operators 
$\tilde\Am_i(\muv)$, $\tilde \Bm(\muv)$
via \eqref{eq:intrusive ops1} and \eqref{eq:intrusive ops2}. In practice, this access is not possible in many FOM solvers. Instead, operator inference aims to learn approximations 
$\hat \Am_i(\muv) \revone{\;\in\mbR^{n\times n_i}}$, $\hat \Bm(\muv) \revone{\;\in\mbR^{n\times \Nu}}$
of the intrusive ROM operators 
from data \cite{PEHERSTORFER2016data}.
This data consists of triples of a ROM state snapshot $\bx_k(\muv)\in\mbR^n$, a ROM state time derivative snapshot $\dot{\bx}_k(\muv)\in\mbR^n$ and an input signal snapshot $\uv_k(\muv)$, $k=0,\dots,K-1$. For these triples, 
operator inference minimizes 
the residual of the ROM ODE \eqref{eq:rom efficient} to infer 
the aggregated ROM operator,
\begin{align}
    \min_{\hat \Om(\muv) = \left[
    [\hat\Am_i]_{i\in\inds}
    \;\; \hat \Bm(\muv) \right]} 
    \sum_{k=0}^{K-1} \left\|
    \hat \Om(\muj) \pv(\bx_k(\muj),\uv_k(\muj))
    -\dot{\bx}_k(\muj)
    \right\|_2^2.
    \label{eq:opinf minimi vecs}
\end{align}
If 
snapshot triples
are 
obtained from the intrusive ROM ODE \eqref{eq:rom efficient} such that \R{4:consistent triple 1}\revfour{$\bx_k(\muj) = \tx(\tk;\muj)$, $\dot\bx_k(\muj) = \dot\tx(\tk;\muj)$} and $\uv_k(\muj) =\uv(\tk;\muj)$ at time steps $\tk$, $k=0,\dots,K-1$,
then 
the intrusive 
aggregated ROM operator
$\tilde \Om(\muv) = \left[
    [\tilde\Am_i]_{i\in\inds}
    \;\; \tilde \Bm(\muv) \right]$ 
is a solution to this minimization problem.

In practice, however, snapshots of the ROM ODE \eqref{eq:rom efficient} are usually not available, but only snapshots of time discretizations of the FOM \eqref{eq:peher fom time conti}, $\x_k(\muj)\in\mbR^N,\; \uv_k(\muj)\in\mbR^\Nu$ corresponding to time steps $\tk$, $k=0,\dots,K$. These snapshots can be used to approximate
the ROM snapshots $\bx_k(\muv)$ in \eqref{eq:opinf minimi vecs} by projecting onto the ROM basis,
\begin{align}
    \bx_k(\muj) = \Vn^T \x_k(\muj), \h k=0,\dots,K,
    \label{eq:proj fom snapshots}
\end{align}
and to approximate the time derivative snapshots $\dot\bx_k(\muv)$ via some finite difference, in the simplest case 
explicit
Euler,
\begin{align}
    \dot \bx_k(\muj) = \frac{\bx_{k+1}(\muj)-\bx_k(\muj)}{t_{k+1}-\tk} , \h k=0,\dots,K-1.
    \label{eq:standard opinf time derivative snapshot}
\end{align}
When
snapshots \eqref{eq:proj fom snapshots} and \eqref{eq:standard opinf time derivative snapshot} 
are used
in the operator inference minimization problem \eqref{eq:opinf minimi vecs}, the inferred 
operator $\hat \Om(\muv)$
 resulting from solving this minimization problem is generally not equal to the intrusive operator $\tilde \Om(\muv)$, but only an approximation.
This approach is \nonintrusive\ because it does not require access to the FOM operators 
$\tilde \Am_i(\muv)$, $\tilde \Bm(\muv)$,
but only the knowledge
about $\inds$ and $\Nu$,
and FOM snapshot data.

Note that operator inference for ROMs with POD bases take two sets of snapshots as input. The first set, the POD snapshot data, consists of FOM state snapshots $\x_k(\muv)\in\mbR^N, \, k = 0,\dots,K_\text{POD}-1$ that are used to compute the POD basis $\Vn$, exactly in the same way as for the corresponding intrusive ROM with POD basis.
The second set, the operator inference snapshot data, consists of triples of a ROM state snapshot $\bx_k\in\mbR^n$, a ROM state time derivative snapshot $\dot{\bx}_k\in\mbR^n$ and an input signal snapshot $\uv_k\in\mbR^\Nu$. This set defines the minimization problem \eqref{eq:opinf minimi vecs}.

Note further that in general the inference of the ROM operators by solving \eqref{eq:opinf minimi vecs} must be performed separately for each choice of the parameter $\muv$, just like the intrusive operator precomputations \eqref{eq:intrusive ops1} and \eqref{eq:intrusive ops2}. For parameter values for which inference data is not available, the corresponding operators are then obtained by interpolation of inferred operators \cite{PEHERSTORFER2016data}.
Under certain assumptions, the parameter dependency of the ROM operators can explicitly be incorporated in the operator inference minimization problem \eqref{eq:opinf minimi vecs} \cite{mcquarrie2023nonintrusive}.
In this work, we
consider the parameter $\muv$ fixed and will omit it in the following sections.

\subsection{Insufficiently informative snapshot data}
\label{sec:insufficient snapshot data}
As has been observed by several researchers \cite{goyal2023guaranteed,mcquarrie2021data,peherstorfer2020sampling}, standard operator inference via \eqref{eq:opinf minimi vecs} often fails to accurately reproduce the intrusive ROM operators. 
There are several causes of such failures such as inconsistency between the FOM time discretization and the discretization of the time derivative snapshots $\dot \bx_k(\muv)$ or non-Markovian information in the snapshot data \cite{peherstorfer2020sampling}.
In this work, we want to highlight 
the remaining cause that has not been resolved yet:
insufficiently informative snapshot data. For this purpose, we write \eqref{eq:opinf minimi vecs} equivalently as
\begin{align}
    \min_{\hat \Om\in\mbR^{n\times \nf}} \left\|\hat \Om \Pm - \dot\bX\right\|_F^2,
    \label{eq:opinf minimi matrices}
\end{align}
with
\begin{align}
    \Pm := 
    \left[ 
    \pv(\bx_0,\uv_0) \;\;\dots\;\;
    \pv(\bx_{K-1},\uv_{K-1})
    \right]
    \; \in\mbR^{\nf\times K}
    \h \text{ and } \h 
    \dot\bX := \left[
    \dot\bx_0 \;\;\dots\;\; \dot\bx_{K-1}
    \right] \; \in\mbR^{n\times K}.
    \label{eq:def p matrix and dot breve x matrix}
\end{align}
This optimization problem has a unique solution $\hat\Om$ if and only if the matrix $\Pm$ has full rank $\nf$ \cite{peherstorfer2020sampling}. In general, there is no guarantee that this requirement is satisfied, even if the number of snapshots $K$ is much larger than $\nf$.
As described in \cite{peherstorfer2020sampling}, strategies to mitigate such rank deficiency include regularization
and reduction of the ROM dimension $n$. However, these strategies do not guarantee accurate reconstruction of the intrusive ROM operators. Regularization can limit the accuracy of the resulting operators. Indeed, many works on regularized operator inference report significant differences between the \nonintrusive\ ROMs and their intrusive counterparts \cite{sawant2023physics,mcquarrie2023nonintrusive,benner2020operator}. 
Reducing the ROM dimension, on the other hand, can result in a least squares matrix of even smaller rank than the initial matrix. Therefore, iterative reduction of the ROM dimension can result in a 
significantly smaller
dimension $n$ corresponding to an intrusive ROM of insufficient accuracy.

\R{2:a posteriori}In the numerical experiments in \cite{peherstorfer2020sampling}, full rank has been achieved by generating and concatenating multiple trajectories of snapshots with randomized input signals. However, the number and length of the trajectories and the distributions of the random input signals are hyperparameters that have to be tuned, and whether the resulting least squares matrix has full rank can only be checked a posteriori.

%% file: sections/existence_and_uniqueness.tex
\section{Proof of full rank of matrix $\Pm$}
\label{sec:full rank proof}

We prove that the  
ROM states
and input signals proposed in Section \ref{sec:basis of the feature space} guarantee that the matrix $\Pm$ has full rank. For this purpose, we first note that $\Pm$ is composed of the submatrices $\Pm\in\mbR^{\np\times\np}$, $\mat Z\in\mbR^{\np\times \Nu}$, $\mat 0\in\mbR^{\Nu\times\np}$ and $\Um\in\mbR^{\Nu\times \Nu}$,
\begin{align}
    \Pm 
        =
    \left[ 
        \pv(\sx_1,\vt 0) \mdots \pv(\sx_\np,\vt 0) 
        \;\;
        \pv(\vt 0,\uv_1) \mdots \pv(\vt 0,\uv_p)
    \right]
    = \begin{bmatrix}
        \bar \Pm & \mat Z \\ \mat 0 & \Um
    \end{bmatrix}
    .
    \label{eq:def specific p matrix}
\end{align}
The matrix
$\bar \Pm\in\mbR^{\np\times\np}$ is the interpolation matrix corresponding to the multivariate interpolation problem of finding the coefficients $\xcoff_\alpha$ of the polynomial
\begin{align}
    \poly: \mbR^n \rightarrow \mbR, 
    \h \tx \mapsto \sum_{|\alpha|\in\inds} \xcoff_\alpha \tx^\alpha
\end{align}
such that 
\begin{align}
    \poly(\tx_i) = b_i, \h i=1,\dots,\np
\end{align}
for the interpolation nodes $\tx_i\in\cup_{i\in\inds} \setXi$ and arbitrary given interpolation values $b_i\in\mbR$.
All entries of $\mat Z$ are equal to $\poly(\vt 0)$.
$\Um\in\mbR^{\Nu\times \Nu}$ is the $\Nu$-dimensional identity matrix.

Since $\Pm$ is a block-triangular matrix, it has full rank if its diagonal blocks $\bar\Pm$ and $\Um$ have full rank. While the identity matrix $\Um$ trivially has full rank, proving that the polynomial interpolation matrix $\bar\Pm$ has full rank is more involved.

Most works on multivariate polynomial interpolation consider the special case $\inds=\mbN_l := \{0,1,\dots,l\}$ for a given $l\in\mbN$ \cite{gasca1990multivariate,gasca2000polynomial}. In this work, however, we allow \textit{gaps} in the \degreeset\ $\inds$, i.e., not all elements of $\mbN_l$ are necessarily contained in $\inds$. Therefore, we define for
all $n\in\mbN^+$, $l\in\mbN_0$ and any $\inds\subset\mbN_l$,
the set of \textit{gappy polynomials} in $\inds$,
\begin{align}
    \mathcal P^n_\inds :=\left\{
     \poly: \mbR^n \rightarrow \mbR, 
    \h \tx \mapsto \sum_{|\alphav|\in\inds} \xcoff_\alphav \tx^\alphav \h \Bigg| \h c_\alphav\in\mbR, \alphav\in\{0,1,\dots,l\}^n       
    \right\}
\end{align}
Then, the following theorem implies that $\bar \Pm$ has full rank.

\begin{mytheo}[Gappy multivariate polynomial interpolation]
\label{theo:multivariate gappy poly interpol}
     For any $n\in\mbN^+$, $\ls\in\mbN_0$, for any \degreeset\ $\inds\subset\mbN_l$,
     there exists a unique polynomial $p\in\mathcal P^n_ \inds$
     such that
     \begin{align}
         \poly(\qv_j) = b_j, \h j=1,\dots,\np,
         \label{eq:gappy multi poly interpol prop}
     \end{align}
     for the interpolation nodes $\qv_j$ in $\cup_{i\in\inds} \setXi$ and arbitrary interpolation values $b_j$.
\end{mytheo}

\subsection{Outline of proof}
\label{sec:outline of proof}

The special case of \degreeset s without gaps, (i.e.\ $\inds = \mbN_l$), is covered by
an existing result on existence and uniqueness of an interpolating multivariate polynomial \cite{nicolaides1972class}
described
in Section \ref{sec:poly interpol max degree}.
In Section \ref{sec:homogen interpol},
we use this result to prove Theorem \ref{theo:multivariate gappy poly interpol} for the special case that $\inds$ contains only one element. 
I.e., we can find for all $i\in\inds$ a homogeneous polynomial $p_i\in\mathcal P^n_{\{i\}}$ that interpolates in the nodes in $\setXi$. 
We use this result as base case in a proof by induction in Section \ref{sec:multivariate gappy poly interpol} to prove the general case where $\inds$ contains multiple elements.
There, the induction hypothesis states that Theorem \ref{theo:multivariate gappy poly interpol} holds for a fixed \degreeset\ $\indsminus$, i.e., there exists a unique polynomial $\pminus$ that interpolates in the interpolation nodes in $\cup_{i\in\indsminus} \setXi$. In the induction step, we show that Theorem \ref{theo:multivariate gappy poly interpol} also holds for the new \degreeset\ $\indsplus:= \{\inew\} \cup \indsminus$ for any $\inew\in\mbN_0$ with $\inew<\min \indsminus$. For this purpose, we add a new polynomial $\pplus$ that has two key properties. First, $\pplus$ attains suitable values in the new interpolation nodes in $\setXx{\inew}$ such that the resulting polynomial $\pminus+\pplus$ attains the prescribed interpolation values $b_j$. Second, $\pplus$ attains 0 in the old interpolation nodes in $\cup_{i\in\indsminus} \setXi$ to not interfere with the interpolation property of $\pminus$ in these nodes. In Section \ref{sec:univariate gappy poly interpol}, we show that such a polynomial $\pplus$ exists in $\mathcal P^n_ {\indsplus}$.

\subsection{Multivariate polynomial interpolation 
without gaps}
\label{sec:poly interpol max degree}

The following theorem is given in \cite{nicolaides1972class} for the special case without gaps in the \degreeset, so $\inds=\mbN_l$.

\begin{mytheo}[Multivariate polynomial interpolation 
without gaps]
\label{theo:interpol with max degree}
    \R{1:Xv 1}
    For any $\n\in\mbN_0$, let $\Xv_0,\Xv_1,\dots,\Xv_\n\in\mbR^\n$ be the vertices of a non-degenerate $\n$-dimensional simplex 
    in $\mbR^\n$. For any $\ls\in\mbN_0$, we define 
    the lattice
\begin{align}
    B(l,\n) :=\left\{
    \xv\in\mbR^\n \Bigg| \xv = \sum_{i=0}^{\n} \lambda_i \Xv_i, \; \lambda_i\in\mbN_0, \; \sum_{i=0}^{\n} \lambda_i = \ls
    \right\}.
    \label{eq:def lattice B l m}
\end{align}
Then there exists a unique polynomial $\poly\in\mathcal P^\n_{\mbN_l}$
such that
\begin{align}
    \poly(\qv_i) = \ival_i, \h i = 1,2,\dots,{\n+l\choose \n}
    ,
\end{align}
for the interpolation nodes $\qv_i$ in  $B(l,\n)$ and arbitrary interpolation values $\ival_i\in\mbR$.
\end{mytheo}

\subsection{Homogeneous multivariate polynomial interpolation}
\label{sec:homogen interpol}

Next we consider the special case that the \degreeset\ $\inds$ contains only one element $l$, so the interpolant is a homogeneous polynomial in $\mathcal P^n_{\{\ls\}}$ and the interpolation nodes are given by $\setXx \ls$. 

\begin{mytheo}[Homogeneous multivariate polynomial interpolation]
\label{theo:homogen poly}
For any $n\in\mbN^+$, $\ls\in\mbN_0$, 
there exists a unique polynomial $\poly\in\mathcal P^n_{\{\ls\}}$ 
such that
\R{1 n_l 1}
\begin{align}
    \poly(\qv_i) = \ival_i, \h i = 1,2,\dots,
    \revone{
    {n+l-1\choose l},}
    \label{eq:interpol homogeneous}
\end{align}
for the interpolation nodes $\qv_i$ in 
$\setXx{l}$ as defined in \eqref{eq:def set X i} and arbitrary interpolation values $\ival_i\in\mbR$.
\end{mytheo}

To prove this theorem, we trace it back to  
Theorem \ref{theo:interpol with max degree} on multivariate polynomial interpolation without gaps.
First, we note that $\setXx l$ lies in an ($n-1$)-dimensional subspace of $\mbR^n$ and constitutes
a lattice of degree $l$ and dimension $n-1$. Hence, we can construct a multivariate polynomial $p\in\mathcal P^n_{\mbN_{l}}$ that interpolates arbitrary values in $\setXx l$. Then, we use the linear dependence of the nodes in $\setXx l$ to reformulate $p$ into an equivalent polynomial in $\mathcal P^n_{\{l\}}$.

\begin{proof}[Proof of Theorem \ref{theo:homogen poly}]
    \R{1:Xv 2}
    We employ Theorem \ref{theo:interpol with max degree} with $\n=n-1$ and the vertices $\Xv_0 = \vt 0\in\mbR^{n-1}$ and $\Xv_j$ the unit vectors in $\mbR^{n-1}$ for $j=1,\dots,n-1$.
    The theorem states that there exists a polynomial
    \begin{align}
        \polyi: \mbR^{n-1}\rightarrow\mbR, \h \polyi(\yv) = \sum_{|\betav|\leq l} \hat a_\betav \yv^\betav,
    \end{align}
    with $\hat a_\betav\in\mbR$ for all multi-exponents $\betav\in\{0,1,\dots,l\}^{n-1}$ with $|\betav|\leq l$, such that
    \begin{align}
        \polyi(\hatqv_i) = 
        f_i, \h i =1,2,\dots,
        {n-1+l\choose n-1}
    \end{align}
    for the prescribed interpolation values $f_i$ and the interpolation nodes $\hatqv_i$ in $B(l,n-1)$.

    \R{1:Xv 3}
    We define the 
    linear map $\map$ that maps the vertices $\Xv_j$ to the unit vectors in $\mbR^n$. Concretely, $\map(\Xv_j)=\onehot_j$, $j=1,\dots,n-1$ and $\map(\Xv_0) = \onehot_n$.
    In view of the definition of $\setXi$ \eqref{eq:def set X i}, we can rewrite the definition of $B(l,m)$ \eqref{eq:def lattice B l m},
    \begin{align}
        B(l,m) = \left\{
            \sum_{j=1}^l \vt a_j \,\Bigg| \vt a_1,\dots,\vt a_l\in\left\{\Xv_0,\dots,\Xv_m\right\} 
        \right\}.
        \label{eq:other def B}
    \end{align}
    Hence,
    we find
    \begin{align}
        \map(B(l,n-1)) = \left\{
            \sum_{j=1}^l \vt a_j \,\Bigg| \vt a_1,\dots,\vt a_l\in\left\{\map(\Xv_0),\dots,\map(\Xv_{n-1})\right\} 
        \right\} = \setXx l.
        \label{eq:mapped points}
    \end{align}
    Consequently,
    the polynomial
    \begin{align}
        \poly(\xv) := \polyi(\map^{-1}(\xv))
    \end{align}
    satisfies \eqref{eq:interpol homogeneous} with the $\qv_i$ in $\setXx l$. 

    To show that this polynomial $\poly$ is an element of $\mathcal P^n_{\{\ls\}}$, we note that the map $\map$ is given by
    \begin{align}
        \map: \mbR^{n-1} \rightarrow \mbR^n, \; \yv \mapsto \begin{bmatrix}
            I_{n-1} \\ \vt 0^T
        \end{bmatrix} \yv + \begin{bmatrix}
            \vt 0 \\ \ls
        \end{bmatrix},
    \end{align}
    with the inverse map
    \begin{align}
        \map^{-1} : \mbR^n \rightarrow \mbR^{n-1}, \; \xv \mapsto [I_{n-1}\;\; \vt 0] \xv.
    \end{align}
    We apply the definition of $\map^{-1}$ to find
    \begin{align}
        \poly(\xv) = \polyi(\map^{-1}(\xv)) = \sum_{|\betav|\leq \ls} \hat a_\betav (\map^{-1}(\xv))^\betav 
        = \sum_{|\betav|\leq \ls} \hat a_\betav \prod_{i=1}^{n-1} \xs_i^{\hat\alpha_i}\revone{,}
    \end{align}
    \R{1 x_i}\revone{with the entries $x_i$ of $\xv\in\mbR^n$, $i=1,\dots,n$.}
    By definition, all $\xv\in \setXx l$ satisfy
    \begin{align}
        \sum_{i=1}^n \xs_i = \ls,
    \end{align}
    so we can multiply $1 = \frac 1 \ls \sum_{i=1}^n \xs_i$ to get
    \begin{align}
        \poly(\xv) = \sum_{|\betav|\leq \ls} \hat a_\betav \left(\prod_{i=1}^{n-1} \xs_i^{\hat\alpha_i}\right) \left( \frac 1 \ls \sum_{i=1}^n \xs_i \right)^{\ls - \sum_{i=1}^n \hat\alpha_i}.
    \end{align}
    As a result, each summand is a monomial of degree $\ls$, so we can write
    \begin{align}
            \poly(\xv) = \sum_{|\alphav|=\ls} a_\alphav \xv^\alphav
    \end{align}
    with suitable $a_\alphav\in\mbR$, for all $\alphav\in\{0,1,\dots,l\}^n$ with $|\alphav|=\ls$. Thereby, we have proven existence of a polynomial $\poly\in\mathcal P^n_{\{\ls\}}$ that satisfies
    \eqref{eq:interpol homogeneous} for any interpolation values 
    \R{1 n_l 2}
    $\ival_i$, $i=1,2,\dots,
    \revone{{n+l-1\choose l}}
    $. Since the corresponding interpolation matrix is square, existence of a solution implies also uniqueness.
\end{proof}

\subsection{Specific gappy multivariate polynomial interpolation}
\label{sec:univariate gappy poly interpol}

As explained in Section \ref{sec:outline of proof}, the proof of Theorem \ref{theo:multivariate gappy poly interpol} in Section \ref{sec:multivariate gappy poly interpol} involves a polynomial $\pplus\in\mathcal P^n_\indsplus$ with $\indsplus := \{\inew\}\cup \indsminus$ such that 
\begin{align}
         \pplus(\qv_j) &= 0, \h \text{ for all } \qv_j\in
         \mathcal X^i \text{ for all } i\in\indsminus
         \label{eq:pplus interpol x is old}
         \\
         \pplus(\qv_j) &= f_j, \h \text{ for all }
         \qv_j\in\mathcal X^\inew,
         \label{eq:pplus interpol x inew}
\end{align}
for prescribed interpolation values $f_j\in\mbR$. To show that such a polynomial $\pplus$ exists, we make the Ansatz
\begin{align}
    \pplus(\qv) = p_\inew(\qv) \hat p\left(\sum_{k=1}^n q_k\right) \h \forall \qv\in\mbR^n,
    \label{eq:induction step polynomial ansatz}
\end{align}
with $p_\inew\in\mathcal P^n_\inew$ that satisfies \eqref{eq:pplus interpol x inew} according to Theorem \ref{theo:homogen poly}, and a univariate polynomial $\hat p$ that solves the specific interpolation problem with interpolation values 0 and 1 described in the following lemma.


\begin{lemma}[Specific gappy univariate polynomial interpolation]
\label{lem:spec gappy uni}
For any $l\in\mbN_0$, for any $\indsplus:=\{\inew\}\,\cup\,\indsminus\subset\mbN_l$ with $\inew<\min\indsminus$, there  exists a 
polynomial $p\in\mathcal P^1_{\indsmi}$
with $\indsmi:= \{i-\inew\big|i\in\indsplus\}$
such that
     \begin{align}
         \poly(i) &= 0, \h \text{ for all }i\in\indsminus \text{ and }\\
         \poly(\inew) &= 1.
     \end{align}    
\end{lemma}

To prove this lemma, we show in Lemma \ref{lem:spec to gen} that for fixed $l\in\mbN_0$, the existence of solutions to such specific interpolation problems with interpolation values 0 and 1 as described above implies the existence of solutions to the general interpolation problems with arbitrary interpolation values. We use this result to prove Lemma \ref{lem:spec gappy uni} with a case discrimination: In Lemma \ref{lem:spec for arbi l}, we consider the case $\inew\neq 0$ and in Lemma \ref{lem:spec gappy interpol with 0} the case $\inew=0$.

\begin{lemma}[Specific and general gappy univariate interpolation]
\label{lem:spec to gen}
If for
some $\ls\in\mbN_0$, 
for any $\Is\subset\Nl
:= \{0,1,\dots,\ls\}$ 
and any pairwise distinct nodes $x_1,\dots,x_\nI\in\mbR$, 
there exists 
a polynomial $\poly\in\PI$ such that
\begin{align}
    \poly(x_j) &= 0 \h \text{ for } j=1,\dots,\nI-1, \text{ and } 
    \label{eq:zero interpol} \\
    \poly(x_\nI) &= 1,
    \label{eq:one interpol}
\end{align}
then there also exists for any $\Is\subset\Nl$ a unique polynomial $\polyii\in\PI$ such that
\begin{align}
    \polyii(x_j) = \ival_j, \h \text{ for } j=1,\dots,\nI,
    \label{eq:general interpol}
\end{align}
for arbitrary interpolation values 
$\ival_j\in\mbR$.
\end{lemma}
\begin{proof}
    Let $i_1<i_2<\dots i_\nI$ be the elements of a given $\Is\subset\Nl$. For $k=1,\dots,\nI$, we define the set $\Is_k$ of the smallest $k$ of these elements. For each of these sets $\Is_k$, let $\poly_k$ be the polynomial in $\mathcal P^1_{\Is_k}$ that satisfies 
    \begin{align}
            \poly_k(x_j) &= 0 \h \text{ for } j=1,\dots,k-1, \text{ and } 
            \label{eq:zero interpol k} \\
            \poly_k(x_k) &= 1.
            \label{eq:one interpol k}
    \end{align}

    To find the polynomial $\polyii\in\PI$ that satisfies \eqref{eq:general interpol}, we make the Ansatz
    \begin{align}
        \polyii(\xs) = \sum_{k=1}^{\nI} \cs_k \poly_k,
    \end{align}
    to obtain an interpolation problem described by the linear system
    \begin{align}
        \intmat \cv = \fv,
        \label{eq:lin sys interpol}
    \end{align}
    with $\cv = [\cs_1 \;\; \dots \;\; \cs_{\nI}]^T \in\mbR^\nI$, $\fv = [\ival_1 \;\;\dots\;\; \ival_\nI]^T\in\mbR^\nI$ and the entries of the matrix $\intmat\in\mbR^{\nI\times\nI}$,
    \begin{align}
        \intmat_{jk} := \poly_k(x_j), \h\h j,k=1,\dots,\nI.
    \end{align}  
    Because of \eqref{eq:zero interpol k} and \eqref{eq:one interpol k}, $\intmat$ is a lower triangular matrix with entries 1 on its diagonal. Hence, $\intmat$ has full rank and the linear system \eqref{eq:lin sys interpol} has a unique solution $\cv$ for any right-hand side $\fv$.
\end{proof}

The following two lemmas extend this result to arbitrary degrees $l$; Lemma \ref{lem:spec for arbi l} for the case that all nodes $x_j$ are positive, Lemma \ref{lem:spec gappy interpol with 0} for the case that $x_\nI=0$ and $0\in\inds$.

\begin{lemma}[Specific gappy univariate interpolation for arbitrary $\ls$]
\label{lem:spec for arbi l}
For any $\ls\in\mbN_0$, for any $\Is\in\Nl:= \{0,1,\dots,\ls\}$, there exists a polynomial $\poly\in\PI$ such that the interpolation properties \eqref{eq:zero interpol} and \eqref{eq:one interpol} hold for any pairwise distinct, positive nodes $x_1,\dots,x_\nI\in\mbR^+$.
\end{lemma}

\begin{proof}
    The proof proceeds by induction along $\ls$.
\begin{itemize}

    \item 
    Base case: For $\ls=0$, $\poly(\xs) = 1$ satisfies the interpolation properties for the only existing set $\Is=\{0\}$.
    
    \item 
    Induction hypothesis: For fixed $\ls$, there exists $\poly\in\PI$ such that the interpolation properties \eqref{eq:zero interpol} and \eqref{eq:one interpol} hold for any $\Is\in\Nl$ and any pairwise distinct, positive nodes $x_1,\dots,x_\nI\in\mbR$.
    
    \item 
    Induction step: We consider a set $\Is\subset\mbN_{\ls+1}:=\{0,1,\dots,\ls,\ls+1\}$ and the pairwise distinct, positive nodes $x_1,\dots,x_\nI\in\mbR^+$. If $\ls+1\notin\Is$, then the induction hypothesis directly applies. Otherwise, if $\ls + 1\in\inds$, then the induction hypothesis holds for $\Is^- := \Is\setminus \{\ls + 1\}$ and the nodes $x_1,\dots,x_{\nI-1}$. Hence, Lemma \ref{lem:spec to gen} implies that there exists a polynomial $\polyii\in\mathcal P^1_{\Is^-}$ such that 
    \begin{align}
    \polyii(x_j) = \ival_j, \h \text{ for } j=1,\dots,\nI-1,
    \end{align}
    for arbitrary interpolation values $\ival_i\in\mbR$.
    Consequently, there exists $\polyii\in\mathcal P^1_{\Is^-}$ such that
    \begin{align}
        \polyii(x_j) = - x_j^{\ls+1}, \h \text{ for } j=1,\dots,\nI-1.
    \end{align}
    Then,
    \begin{align}
        \hat p(x) := \bar p(x) + x^{\ls+1}
    \end{align}
    satisfies \eqref{eq:zero interpol}, i.e.,  $x_1,\dots,x_{\nI-1}$ are roots of $\hat p$. By Descartes' rule of signs, $\hat p$ has at most $\nI-1$ positive roots, hence $\hat p(x_\nI) \neq 0$. Consequently, we can normalize $\hat p$ to find the desired polynomial,
    \begin{align}
        \poly(\xs) := 
        \frac{\hat p(x)}{\hat p(x_\nI)}.
        \label{eq:def specific polynomial}
    \end{align}
\end{itemize}
\end{proof}

\begin{lemma}[Specific gappy univariate interpolation for arbitrary $l$ with $\x_\nI=0$]
\label{lem:spec gappy interpol with 0}
    For any $\ls\in\mbN_0:= \{0,1,\dots,\ls\}$, for any $\Is\in\Nl$ with $0\in\inds$, there exists a polynomial $\poly\in\PI$ such that the interpolation properties \eqref{eq:zero interpol} and \eqref{eq:one interpol} hold for any pairwise distinct, positive nodes $x_1,\dots,x_{\nI-1}\in\mbR^+$ and $x_\nI=0$.
\end{lemma}
\begin{proof}
    Let $\inds^- := \inds \setminus\{0\}$.
    According to Lemma \ref{lem:spec to gen} and Lemma \ref{lem:spec for arbi l}, there exists a polynomial $\bar p\in\mathcal P^1_{\inds^-}$ such that
    \begin{align}
        \bar p(x_j) = -1
        \h \text{ for } j=1,\dots,\nI-1.
    \end{align}
    Then,
    \begin{align}
        p(x) := \bar p(x) + 1
    \end{align}
    is element of $\PI$ and satisfies \eqref{eq:zero interpol}.
    Because $\bar p\in\mathcal P^1_{\inds^-}$ with $0\notin\inds^-$, $x$ is a factor of $\bar p$, so $\bar p(0) = 0$. Consequently $p$ also satisfies \eqref{eq:one interpol}.
\end{proof}

With these preparations, we can prove 
Lemma \ref{lem:spec gappy uni}.

\begin{proof}[Proof of Lemma \ref{lem:spec gappy uni}]
    For any $l\in\mbN$, for any $\hat \inds = \{\inew\}\cup\indsminus\subset\mbN_l$ with $\inew<\min\indsminus$, let $\inds = \{i-\inew\big|i\in\hat\inds\}$ and $\{x_1,\dots,x_\nI\}=\hat \inds$. Then Lemma \ref{lem:spec for arbi l} implies the statement if $\inew\neq 0$ and Lemma \ref{lem:spec gappy interpol with 0} if $\inew=0$. 
\end{proof}

    

\subsection{Gappy multivariate polynomial interpolation}
\label{sec:multivariate gappy poly interpol}

Now we have all ingredients to prove Theorem \ref{theo:multivariate gappy poly interpol}.

\begin{proof}[Proof of Theorem \ref{theo:multivariate gappy poly interpol}]
The proof is performed by induction along the number of elements in $\inds$.
\begin{itemize}
    \item Base case: Let $\inds=\{i\}$ for any $i\in\mbN_l$ with arbitrary $l\in\mbN_0$, then Theorem \ref{theo:homogen poly} implies the existence and uniqueness of a polynomial $p\in\mathcal P^n_{\{i\}}$ that satisfies the interpolation property \eqref{eq:gappy multi poly interpol prop} in the interpolation nodes in $\setXi$.
    \item Induction hypothesis: Let for a fixed \degreeset\ $\indsminus\subset\mbN_l$ with arbitrary $l\in\mbN_0$ exist a unique polynomial $p^-\in\mathcal P^n_\indsminus$ that satisfies the interpolation property \eqref{eq:gappy multi poly interpol prop} in the interpolation nodes $\cup_{i\in\indsminus} \setXi$.
    \item Induction step: We consider the \degreeset\ $\indsplus :=\{\inew\}\cup\indsminus$ with arbitrary $\inew\in\mbN_l$. Without loss of generality, we can assume $\inew<\min\indsminus$. As described in Section \ref{sec:univariate gappy poly interpol}, we can construct a polynomial
    \begin{align}
    \pplus(\qv) = p_\inew(\qv) \hat p\left(\sum_{k=1}^n q_k\right) \h \forall \qv\in\mbR^n,
\end{align}
    with $p_\inew\in\mathcal P^n_\inew$ and $\hat p\in\mathcal P^1_\indsmi$, so $p\in\mathcal P^n_\indsplus$.
    Because of Lemma \ref{lem:spec gappy uni} and $\sum_{k=1}^n q_k = i$ for all $\qv\in\setXi$ for any $i\in\mbN_0$ , we can choose $\hat p$ such that
    \begin{align}
        \hat p\left(\sum_{k=1}^n q_k\right) &= 0, \h \text{ for all } \qv\in
         \mathcal X^i \text{ for all } i\in\indsminus
         \\
         \hat p\left(\sum_{k=1}^n q_k\right) &= 1, \h \text{ for all }
         \qv\in\mathcal X^\inew,
    \end{align}
    and $p_\inew$ such that 
    \begin{align}
        p_\inew(\qv_j) = b_j - \pminus(\qv_j), \h \text{ for all } \qv_j\in\setXx \inew,    
    \end{align}
    with the interpolation values $b_j$ prescribed in \eqref{eq:gappy multi poly interpol prop}.
    Consequently, the polynomial $\pplus$ satisfies
    \begin{align}
         \pplus(\qv_j) &= 0, \h \text{ for all } \qv_j\in
         \mathcal X^i \text{ for all } i\in\indsminus \; \text{ and }
         \\
         \pplus(\qv_j) &= b_j - \pminus(\qv_j), \h \text{ for all }
         \qv_j\in\mathcal X^\inew.
\end{align}
As a result, the polynomial $\pminus + \pplus$ satisfies
\begin{align}
        \pminus(\qv_j) + \pplus(\qv_j) &= b_j + 0 = b_j, \h \text{ for all } \qv_j\in
         \mathcal X^i \text{ for all } i\in\indsminus \; \text{ and }
         \\
         \pminus(\qv_j) = \pplus(\qv_j) &= \pminus(\qv_j) + b_j - \pminus(\qv_j) = b_j, \h \text{ for all }
         \qv_j\in\mathcal X^\inew,
         \end{align}
and hence satisfies the interpolation property \eqref{eq:gappy multi poly interpol prop}. Since the corresponding interpolation matrix is square, existence of a solution also implies uniqueness.
\end{itemize}

\end{proof}